\documentclass[a4paper,12pt]{amsart}
\usepackage[utf8]{inputenc}
\usepackage[english]{babel} 

\usepackage[T1]{fontenc}

\usepackage{amsthm}
\usepackage{amsmath}
\usepackage{amsfonts}
\usepackage{mathtools}
\usepackage[colorlinks]{hyperref}
\usepackage{enumerate}
\usepackage{amssymb}
\usepackage{tikz,pgf}
\usepackage{graphicx}
\usetikzlibrary{matrix,arrows,decorations.pathmorphing}
\usetikzlibrary{calc}
\makeatletter
\newbox\dottedarrow@box
\setbox\dottedarrow@box\hbox
  {%
    \begin{tikzpicture}
      \draw[dotted,->] (0,0) -- (1.5em,0);
    \end{tikzpicture}%
  }
\newcommand*\dottedarrow
  {\relax\ifmmode\expandafter\dottedarrow@m\else\expandafter\dottedarrow@t\fi}
\newcommand*\dottedarrow@t[1][1.5em]
  {\resizebox{#1}{!}{\raisebox{.5ex}{\usebox\dottedarrow@box}}}
\newcommand*\dottedarrow@m[1][]
  {%
    \if\relax\detokenize{#1}\relax
      \mathchoice
        {\dottedarrow@t}
        {\dottedarrow@t}
        {\dottedarrow@t[1.1em]}
        {\dottedarrow@t[0.9em]}%
    \else
      \dottedarrow@t[#1]%
    \fi
  }
\makeatother

\def\bdi{\begin{diagram}}
\def\edi{\end{diagram}}

\newtheorem{thm}{Theorem}[section]
\newtheorem{cor}[thm]{Corollary}
\newtheorem{lem}[thm]{Lemma}
\newtheorem{prop}[thm]{Proposition}
\theoremstyle{definition}
\newtheorem{defi}[thm]{Definition}
\newtheorem{defis}[thm]{Definitions}
\newtheorem{conj}[thm]{Conjecture}
\newtheorem{conv}[thm]{Convention}
\newtheorem{nota}[thm]{Notation}
\newtheorem{rem}[thm]{Remark}
\newtheorem{rems}[thm]{Remarks}
\newtheorem{exa}[thm]{Example}
\newtheorem{exas}[thm]{Examples}
\newtheorem{prob}[thm]{Problem}
\newtheorem{probs}[thm]{Problems}
\newtheorem{ques}[thm]{Question}
\newtheorem{sett}[thm]{Setting}
\newtheorem{sit}[thm]{}

\newcommand{\Spec}{\operatorname{{\rm Spec}}}

\newcommand{\Bl}{ \operatorname{{\rm Bl}}}

\newcommand{\supp}{ \operatorname{{\rm supp}}}

\def\pr{\mathop{\rm pr}}

\renewcommand{\epsilon}{\varepsilon}

\def\and{\quad\mbox{and}\quad}

\newcommand{\F}{\ensuremath{\mathbb{F}}}

\newcommand{\Z}{\ensuremath{\mathbb{Z}}}

\newcommand{\G}{\ensuremath{\mathbb{G}}}

\newcommand{\A}{\ensuremath{\mathbb{A}}}

\newcommand{\bK}{{\ensuremath{\rm \bf K}}}
\newcommand{\bk}{{\ensuremath{\rm \bf k}}}

\newcommand{\hI}{{\hat I}}

\newcommand{\hD}{{\hat D}}

\newcommand{\hG}{{\hat G}}
\newcommand{\hH}{{\hat H}}

\newcommand{\hX}{{\hat X}}

\newcommand{\hZ}{{\hat Z}}

\newcommand{\hE}{{\hat E}}

\newcommand{\hO}{{\hat O}}

\newcommand{\hQ}{{\hat Q}}

\newcommand{\hf}{{\hat f}}
\newcommand{\hh}{{\hat h}}

\newcommand{\hx}{{\hat x}}
\newcommand{\hy}{{\hat y}}

\newcommand{\hv}{{\hat v}}
\newcommand{\hw}{{\hat w}}

\newcommand{\tv}{{\tilde v}}

\newcommand{\tC}{{\tilde C}}
\newcommand{\tD}{{\tilde D}}

\newcommand{\tH}{{\tilde H}}

\newcommand{\tQ}{{\tilde Q}}

\newcommand{\tX}{{\tilde X}}

\newcommand{\tE}{{\tilde E}}
\newcommand{\tP}{{\tilde P}}
\newcommand{\tS}{{\tilde S}}

\newcommand{\tp}{{\tilde p}}
\newcommand{\tq}{{\tilde q}}
\newcommand{\ts}{{\tilde s}}

\newcommand{\tx}{{\tilde x}}
\newcommand{\ty}{{\tilde y}}

\newcommand{\hY}{{\hat Y}}

\newcommand{\brE}{{\breve E}}

\newcommand{\brH}{{\breve H}}

\newcommand{\brP}{{\breve P}}
\newcommand{\brQ}{{\breve Q}}

\newcommand{\brX}{{\breve X}}

\newcommand{\bY}{{\bar Y}}

\newcommand{\ddR}{{{\ddot R}}}
\newcommand{\ddU}{{{\ddot U}}}

\newcommand{\ddc}{{{\ddot c}}}
\newcommand{\ddx}{{{\ddot x}}}
\newcommand{\ddy}{{{\ddot y}}}

\def\fo{{\mathfrak o}}

\newcommand{\cL}{{\ensuremath{\mathcal{L}}}}

\newcommand{\cF}{{\ensuremath{\mathcal{F}}}}

\newcommand{\cS}{{\ensuremath{\mathcal{S}}}}
\newcommand{\cE}{{\ensuremath{\mathcal{E}}}}
\newcommand{\cA}{{\ensuremath{\mathcal{A}}}}
\newcommand{\cJ}{{\ensuremath{\mathcal{J}}}}
\newcommand{\cO}{{\ensuremath{\mathcal{O}}}}

\newcommand{\cD}{{\ensuremath{\mathcal{D}}}}

\newcommand{\cI}{{\ensuremath{\mathcal{I}}}}

\newcommand{\cP}{{\ensuremath{\mathcal{P}}}}

\newcommand{\id}{{\rm id}}

\renewcommand{\rho}{\varrho}

\def\bals#1\eals{\begin{align*}#1\end{align*}}
\def\bal#1\eal{\begin{align}#1\end{align}}

\def\A{{\mathbb A}}

\def\PP{{\mathbb P}}

\renewcommand{\phi}{\varphi}

\newcommand{\bnum}{\begin{enumerate}}
\newcommand{\enum}{\end{enumerate}}

\newcommand{\brem}{\begin{rem}}
\newcommand{\brems}{\begin{rems}}
\newcommand{\erem}{\end{rem}}
\newcommand{\erems}{\end{rems}}
\newcommand{\bprob}{\begin{prob}}
\newcommand{\eprob}{\end{prob}}
\newcommand{\bprobs}{\begin{probs}}
\newcommand{\eprobs}{\end{probs}}
\newcommand{\bques}{\begin{ques}}
\newcommand{\eques}{\end{ques}}
\newcommand{\bexa}{\begin{exa}}
\newcommand{\bexas}{\begin{exas}}
\newcommand{\eexa}{\end{exa}}
\newcommand{\eexas}{\end{exas}}
\newcommand{\bdefi}{\begin{defi}}
\newcommand{\edefi}{\end{defi}}
\newcommand{\bdefis}{\begin{defis}}
\newcommand{\edefis}{\end{defis}}
\newcommand{\bcor}{\begin{cor}}
\newcommand{\ecor}{\end{cor}}
\newcommand{\blem}{\begin{lem}}
\newcommand{\elem}{\end{lem}}
\newcommand{\bconv}{\begin{conv}}
\newcommand{\econv}{\end{conv}}
\newcommand{\bconj}{\begin{conj}}
\newcommand{\econj}{\end{conj}}
\newcommand{\bprop}{\begin{prop}}
\newcommand{\eprop}{\end{prop}}
\newcommand{\bthm}{\begin{thm}}
\newcommand{\ethm}{\end{thm}}
\newcommand{\bnota}{\begin{nota}}
\newcommand{\enota}{\end{nota}}
\newcommand{\bsit}{\begin{sit}}
\newcommand{\esit}{\end{sit}}
\newcommand{\be}{\begin{equation}}
\newcommand{\ee}{\end{equation}}
\newcommand{\bproof}{\begin{proof}}
\newcommand{\eproof}{\end{proof}}
\newcommand{\bsett}{\begin{sett}}
\newcommand{\esett}{\end{sett}}
\def\ba{\begin{array}}
\def\ea{\end{array}}

\newcommand{\nlin}{\unitlength1mm\begin{picture}(0,9.25)
                      \put(0,0.75){\line(0,1){8.5}}
                     \end{picture}}

\newcommand{\vlin}[1]{\hspace{0.75mm}\unitlength1mm\begin{picture}
(#1,0)
                      \put(0,0){\line(1,0){#1}}
                     \end{picture}\hspace{0.75mm}\rule[-3mm]{0mm}
                     {4mm}}

\def\llin{\vlin{11.5}}
\newcommand{\lin}{\vlin{8.5}}

\newcommand{\cou}[2]{\unitlength1mm\begin{picture}(0,8)
   \put(0,0){\circle{1.5}}
   \put(0,3){\makebox(0,5)[b]{$#1$}}
   \put(0,-7){\makebox(0,4)[t]{$#2$}}
     \end{picture}
     \rule[-7mm]{0mm}{7mm}}

\newcommand{\crl}[2]{\unitlength1mm\begin{picture}(0,8)
   \put(0,0){\circle{1.5}}
   \put(-5,0){\makebox(0,5)[b]{$#1$}}
  \put(5,0){\makebox(0,5)[b]{$#2$}}
     \end{picture}
     \rule[-7mm]{0mm}{7mm}}

\newcommand{\cshiftup}[2]{\unitlength1mm\begin{picture}(0,9.25)
                      \put(0,10){\crl{#1}{#2}}
                     \end{picture}}

\begin{document}
\title[Gromov's ellipticity of prinicpal $\G_m$-bundles]{
Gromov's ellipticity of principal $\G_m$-bundles}

\author{Sh. Kaliman}


\thanks{2020 \emph{Mathematics Subject Classification} 
Primary  14J60, 14M25, 14M27; Secondary 32Q56} 
\keywords{ Gromov ellipticity, spray, 
uniformly rational variety, complete surface}
\date{}
\maketitle
\begin{abstract} We prove that every nontrivial  principal $\G_m$-bundle over
a complete stably uniformly rational variety is algebraically elliptic in the sense of Gromov.\end{abstract}
\section{Introduction}\label{sec:intro}
All varieties in this paper are smooth algebraic 
varieties defined over an algebraically closed field  $\bk$ of characteristic zero.
This work is  continuation of 
\cite{AKZ23} and \cite{KZ23b} and, thus, we adhere
to the terminology and notations of these papers (see also \cite{For17} 
for the general theory of Oka manifolds and  \cite{For23} for recent advances in this area). 
In particular, when we talk about Gromov's ellipticity
(or, ellipticity for short) we always mean {\bf algebraic} ellipticity leaving the holomorphic case aside.
This enables us to use the following advantage - the notions of algebraic Gromov's ellipticity
and subellipticity coincide \cite{KZ23a}, while in the holomorphic case it is still an open problem.
Recall that an algebraic variety is called
uniformly rational if every point in it has a Zariski open neighborhood isomorphic to an open subset of $\A^n$
and it is still unknown whether a smooth 
complete rational variety is uniformly rational (see, \cite[3.5.E$'''$]{Gro89},
\cite[Question 1.1]{BB14} and \cite[p. 41]{CPPZ21}). 
Recall also that a variety $X$ is stably uniformly rational if $X\times \A^m$ is uniformly rational for some $m\geq 0$.
The main results of
\cite[Theorem 3.3, Corollary 3.7]{AKZ23} are the following.

\bthm\label{int.t1} Every complete stably uniformly rational variety is Gromov's elliptic.
\ethm

\bthm\label{int.t2} Let $X$ be a complete uniformly 
rational variety and $\cD \to X$ be an ample or 
anti-ample line bundle on $X$
with zero section $Z_\cD$.
Then $\cD\setminus Z_\cD$ is elliptic. 
\ethm

However, the assumption that $\cD\to X$ is ample or anti-ample is not necessary (e.g., see \cite[Theorem 3.8]{AKZ23}).
What is necessary, of course, is the requirement that $\cD\to X$ is a nontrivial line bundle.
This leads to the following.\\

{\bf Question 1.} {\em Let $X$ be a complete stably uniformly 
rational variety and $\cD\to X$ be a nontrivial  line bundle on $X$
with zero section $Z_\cD$. Is $Y=\cD\setminus Z_\cD$ elliptic?}\\

In this paper we show that the answer to this question is positive (Corollary \ref{main.c1}).  Let us briefly discuss
our approach in the crucial case when $X$ is a smooth complete rational surface (so 
$X$ is Gromov's elliptic by   Theorem \ref{int.t1}).
Given $v \in X$ we look for   an open set $B\subset \A^1$  and a birational morphism $\phi : B\times \PP^1 \to X$  locally invertible  at $v$. In particular, $v$ is contained in $C =\phi (b_0\times \PP^1)$ for some $b_0\in B$.
 It follows from \cite[Corollary 3.8]{KZ23b} 
and  \cite[Corollary 2.7]{AKZ23} that
 $Y$ admits a family of rank 1 sprays dominating  on the fiber over $v$ if the restriction of $\cD$ yields a nontrivial 
 line bundle on $C$. If this is true for every $v \in X$, then one has subellipticity of $Y$. 
 Thus, since subellipticity implies ellipticity \cite{KZ23a} we only need to check that
 $\cD|_C$ is not trivial. Letting $\cD=\cO_X(D)$ for some SNC divisor $D$ we observe that
 $\cD|_C$ is not trivial if $C$ meets only one irreducible component $H$ of $D$ (where $C\ne H$). Hence, our aim
 is to find such $\phi, D$ and $H$ for any given $v \in X$.
 The realization of this plan is  executed in Sections \ref{con}-\ref{it} (after preliminaries
 in Section \ref{pre})  which concludes the case of surfaces. 
 
 In Section \ref{main} we extend this result
 to the case of higher dimensions using the following scheme. 
 By a trick of L\'arusson (see  Proposition \ref{main.p3})
  it suffices to consider the case of uniformly rational $X$ only.
 For a nontrivial line bundle $\cD =\cO_X (D)$ and every $v \in X$
 we need again to find a smooth variety $V$
of dimension $n-1$ and
 a morphism $\phi : V\times \PP^1 \to X$ invertible over
 a neighborhood $U$ of  $v$ such that the image of $\phi$ meets only one irreducible component of $D$.
Such $\phi$ can be easily found if there is a birational morphism $X\to\PP^n$ \'etale on $U$
for which the pushforward of $D$ generates a nontrivial line bundle on $\PP^n$.
Otherwise, for any morphism $\kappa : \hX \to X$ from a smooth variety $\hX$ invertible over $U$
  we observe that it suffices to prove the existence of  $\phi$
   for the triple $(\hX, \hat \cD, \hv)$ instead of $(X,\cD, v)$
  where $\hv=\kappa^{-1} (v)$ and $\hat \cD = \kappa^*\cD$ coincides with $\cO_\hX (\hD)$. 
  Then using the induction  on dimension
  we show that for an appropriate $\hX$ there is morphism $\gamma : \hX \to Z\simeq B\times Z_0 $
where $Z_0$ is a smooth projective rational surface 
 such that the following is true (Proposition \ref{main.p1} and Lemma \ref{main.l8}):

 (a) there is a codimension 2 subvariety $M\subset Z$ for which $\gamma$ is invertible over $Z\setminus M$
 and $M$ meets every surface $b \times Z_0, \, b \in B$ at  most at a finite set $M_b$;
 
 (b) $\gamma (\hD)$ is a disjoint union of a subset of $M$ and a hypersurface of the form $B\times E$
 where $E$ is an SNC curve in $Z_0$.
 
For some $b_0\in B$ we have $\gamma (\hv) \in b_0\times Z_0\simeq Z_0$.
Using the case $\dim X=2$ one can find a morphism
 $\psi_0 : B_1\times \PP^1 \to Z_0$ invertible over a neighborhood of $\gamma (\hv)$
 such that the image of $\psi_0$ meets $E$ only at one irreducible
 component $E_0$ and it does not meet $M_{b_0}$ (Proposition \ref{it.p3}). Letting $b_1^0\in B_1$ be such that
 $\gamma (\hv) \in \psi_0 (b_1^0\times \PP^1)$ we see that there is a neighborhood $W$
 of $(b_0,b_1^0)$ in $B\times B_1$ for which the image of the morphism $W\times \PP^1 \to Z,
 (b,b_1, p) \mapsto (b, \psi_0 (b_1,p))$ meets $B\times E$ at the component  $B\times E_0$ only and
 it does not meet $M$.
 Then (a) allows a lift of this morphism to a desired morphism $W\times \PP^1 \to \hX$.
 
\section{Preliminaries}\label{pre}

The following basic definitions can be found in \cite{Gro89} and \cite[Definition 6.1.1]{For17}.

\bdefi\label{pre.d1}
 A \emph{spray of rank $r$} over a smooth algebraic variety $X$ is a triple $(E,p,s)$ consisting of a vector bundle 
$p\colon E\to X$ of rank $r$ and a morphism $s\colon E\to X$ such that $s|_Z=p|_Z$ where $Z\subset E$ stands 
for the zero section of $p$. This spray is \emph{dominating at $x\in X$} if the restriction $s|_{E_x}\colon E_x\to X$ 
to the fiber $E_x=p^{-1}(x)$ is dominant at the origin $0_x$ of the vector space $E_x$.
The variety $X$ is called \emph{elliptic} if it admits a spray $(E,p,s)$ which is dominating at each point $x\in X$.
The variety $X$ is called 
\emph{subelliptic} if it
admits a family of sprays $(E_i,p_i,s_i)$ defined over $X$ which is dominating at each point $x\in X$, that is,
\[T_xX=\sum_{i=1}^n {\rm d}s_i (T_{0_{i,x}} E_{i,x}) \quad \forall x\in X.\]
\edefi

Recall  the following \cite[Definition 2.5]{AKZ23} (see, also \cite[Definition 2.7]{KZ23b}). 

\bdefi\label{pre.d2}
We say that a 
complete rational curve $C$ on a smooth variety $X$
verifies the  \emph{strengthened two-orbit property} 
at a smooth point $v\in C$  
if  there exists a pair  of rank 1 
sprays $(E_i,p_i,s_i)$ $(i=1,2)$ 
on $X$ such that $C$ is covered by the 
one-dimensional $s_i$-orbits 
$O_{i,v}$,  where $s_i\colon p_i^{-1}(v)\to  O_{i,v}$ 
is a birational morphism \'etale over  $v$
 and $(E_i,p_i,s_i)$ restricts to a spray 
on $O_{i,v}$ dominating at~$v$.

If for any $v\in X$ there exists a curve 
$C=C_v$ as above, 
then we say that $X$ verifies the 
\emph{strengthened curve-orbit property}.
\edefi

The reason  why we need the curve-orbit property is the following \cite[Corollary 2.7]{AKZ23}.

\bprop\label{pre.p1} Let $X$ be a smooth
variety, $\rho : \cD\to X$ be a 
line bundle and $Y=\cD\setminus Z_\cD$. Suppose that
\begin{itemize}
\item[{\rm (i)}] $X$ is elliptic and verifies the strengthened curve-orbit 
property with respect to
a collection $\cF$ of projective  rational curves 
$\{C_v\}_{v\in X}$ on $X$ that are smooth at $v$,
 and
\item[{\rm (ii)}] $\rho\colon \cD \to X$ restricts to a nontrivial line bundle 
on each member $C_v$ of $\cF$.
\end{itemize}
Then $Y$ is elliptic.
\eprop
The first step in the proof of Proposition \ref{pre.p1} is to establish subellipticity using
the strengthened curve-orbit property and then ellipticity follows by \cite{KZ23a}.
\bprop\label{pre.p2} 
Let $X$ be a smooth variety of dimension $n$ and $B$ be a smooth variety of dimension $n-1$.
Assume that $X$ admits a birational morphism
$\phi\colon B\times\PP^1\to X$ locally invertible over a neighborhood of 
a point $v\in X$.
Let $C=\phi (b\times \PP^1), \, b\in B$
be the  curve containing $v$. Then
$C$ verifies the strengthened two-orbit property at $v$. 
\eprop
\bproof Let $v=\phi (b,u_0)$ where $u_0\in \PP^1$. By \cite[Lemma 2.4]{AKZ23} for every point 
$ u\in \PP^1\setminus \{u_0\}$
there is a rank 1 spray $(E_u,p_u, s_u)$ on $X$ whose restriction to $C$ is also a spray such that
$s_u|_{p_u^{-1} (v)}: p_u^{-1} (v)\to C$ is  a birational morphism \'etale over $v$ 
and the $s_u$-orbit is equal to $\phi (b\times (\PP^1\setminus \{u\}))$. This yields
the desired conclusion.
\eproof
Propositions \ref{pre.p1} and \ref{pre.p2} imply now the following.
\bcor\label{pre.c1} Let $X$ be a smooth complete  elliptic variety,
$\rho : \cD\to X$ be a line bundle with the zero section $Z$, and $Y=\cD\setminus Z_\cD$. 
Suppose that for every $v \in X$ there is a morphism $\phi_v : B_v\times \PP^1\to X$ as in Proposition
\ref{pre.p2} and $C_v$ be the curve $\phi_v(b_v\times \PP^1)$ containing $v$.
Let the restriction of $\cD$ to every curve $C_v,\, v\in X$ be a nontrivial line bundle. Then $Y$ is elliptic.
\ecor

\bexa\label{pre.exa1} Let $X, \cD$ and $Y$ be as in Corollary \ref{pre.c1}
with $\dim X=1$ and $\cD$ being nontrivial.
In this case $X\simeq \PP^1$ since any other smooth complete curve is not elliptic. Choosing $\phi : \PP^1 \to X$
to be the identity map we see that $Y$ is elliptic.
\eexa

\section{Contractible graphs}\label{con}

 Let us recall first  terminology and  some standard facts about dual weighted graphs of curves.
Let $G$ be a simple normal crossing (SNC) curve in a smooth  surface $X$,
i.e., all irreducible components
of $G$ are smooth and each of its singularities is a point where exactly two components meet transversely.
The {\em dual graph} of $G$   is the graph whose vertices
are the irreducible components of $G$ and
the edges between the vertices are the intersection points of these components.
Assuming that $G$ is complete one can consider the
 {\em  weighted dual graph} $\Gamma$ of $G$  that is the dual graph of $G$ with each vertex $C$
equipped with the weight equal to the selfintersection number $C \cdot C$
 of this component $C$ in $X$. 
 A vertex $C$ is {\em linear (resp. endpoint)} if it has two (resp. one) neighbors in $\Gamma$.
 If $C$ has at least three neighbors it is called a  {\em  branch point}. A non-circular connected graph without branch
 points will be called {\em linear.}
Suppose that $\pi : \hX \to X$ is the blowing up of a point $p\in G$ and $G' =\pi^{-1} (G)$.
Then $G'$ is again an SNC curve. If $p$ is the intersection of two components of $G$ with weights
$w_1$ and $w_2$,  then the graph of $G'$ is obtained from $\Gamma$ by the following change
$$\lin\cou{}{w_1}\lin\cou{}{w_2} \lin \, \, \, \,\,\, \Longrightarrow \, \, \, \, \, \, \, \lin\cou{}{w_1-1}\lin\cou{}{-1}\llin\cou{}{\, \, \, \, \, \, w_2-1} \lin \, \, \, ,$$ while the rest of the dual graph remains the same (such change is called an {\em inner blowing up}).
If $p$ is a smooth point of $G$ located in a component  with weight $w$, then the dual graph of $G'$
is obtained from $\Gamma$ by the following change\\[2ex]
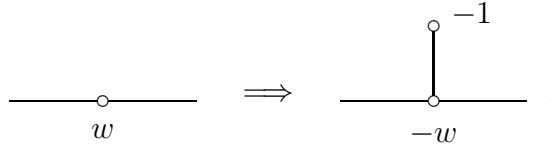
 $$ \llin\cou{}{w}\llin \, \, \, \,\,\, \Longrightarrow  \, \, \, \,\,\,
\llin\cou{\quad \, \, }{-w} \nlin\cshiftup{}{-1}
\llin \, \, \, .$$
The latter change is called an {\em outer blowing up} and unlike the inner blowing up it does not describe
the curve $G'$ uniquely since one needs to indicate the point of $w$ at which the blowing up occurs.
The contraction of a $(-1)$-component $E$ in $G$ also preserves SNC type, provided that $E$
is not a branch point. In this case the image of $G$ has the dual graph obtained from $\Gamma$ by
reversion of the above diagrams. A dual graph is called contractible if such contractions can reduce it
to an empty graph.
\bnota\label{con.n1} From now on by $\A^n_{u_1, \ldots, u_n}$ we denote an affine space $\A^n$
equipped with a coordinate system $(u_1, \ldots, u_n)$.
\enota
\blem\label{con.l1} Let $\Gamma$ be a  contractible linear graph with a unique $(-1)$-vertex $E$
and $C$ be a neighbor of $E$ in $\Gamma$.
Then $\Gamma$ is the dual graph of the minimal resolution $\pi : W\to \A^2_{x,y}$ of the indeterminacy points 
of a rational function $y^m/x^k$ where $k$ and $m\geq 1$ are relatively prime. 
Furthermore,  there are  relatively prime $l$ and $n$ and a neighborhood $U\simeq \A^2_{u,v}$
of $E\cap C$ in $W$ such that $E\cap U$ and $C\cap U$
are coordinate axes and $(u,v)$ coincides with $(y^m/x^k, x^n/y^l)$ or $(x^k/y^m, y^l/x^n)$.
 \elem
\bproof  We use induction on the number of vertices.
Note that every  weighted linear contractible graph with $E$ as an endpoint must be of the form 
$$\cou{C_1}{-1}\lin\cou{C_2}{-2}\lin \ldots \lin\cou{C_k}{-2} \,\,\,.$$
This is exactly the  resolution graph of indeterminacies of $y^m/x^k$ when $m=1$.
Note that each point $C_i\cap C_{i+1}$ has a neighborhood in $W$ isomorphic to $\A^2$ on which
the functions   ${x^{k+1-i}/y}$ and $y/{x^{k-i}}$  yield a coordinate system.
For $k,m>1$ denote now by $C_1$ and $C_2$ the
neighbors of $E$.  Denote their proper transforms in the image $W'$ of $W$ after the contraction of $E$ by
$C_1'$ and $C_2'$ and the new dual graph by $\Gamma'$. Note that one of them (say, $C_1'$) is a $(-1)$-vertex, while the weight of 
$C_2'$ is at most $-2$ to guarantee further contraction. By the induction assumption
one can suppose that $\Gamma'$ is the minimal resolution  $ W'\to \A^2$ of the indeterminacy points of a function
${y^j}/{x^i}$.  Furthermore, there is a neighborhood of the point $C_1'\cap C_2'$ in $W'$ isomorphic to $\A^2$
equipped with the coordinate system  $(u,v) =({y^j}/{x^i}, x^n/{y^l})$ (resp. $(u,v) =(x^i/{y^j}, y^l/x^n)$)
with $C_1'$ given locally by $v=0$ and the matrix
\be \label{con.eq1} \left[ {\begin{array}{cc}
 j& l \\      
i & n
 \end{array} } \right]  \ee
 having the determinant $\pm 1$ (this assumption is obvious in the case $m=1$). 
 The blowing up of $W'$
  at   $C_1'\cap C_2'$ produces $E$ and resolves
 indeterminacy of ${u}/{v} ={y^{j+l}}/{x^{i+n}}$. Note that the assumption on the matrix in Formula \eqref{con.eq1}
 implies that $m=j+l$ and $k=i+n$ are relatively prime and the determinants of the matrices
 \be \label{con.eq2}  \left[ {\begin{array}{cc}
 j+l& l \\      
i+n & n
 \end{array} } \right]    \text{  and  }  \left[ {\begin{array}{cc}
 j+l& j \\      
i+n& i
 \end{array} }  \right]  \ee
 are $\pm 1$. Note also that $(u, v/u)$ (resp. $(v, u/v)$) is a local coordinate system in a neighborhood
 of $C_1\cap E$ (resp. $C_2\cap E$).
This concludes the proof.
\eproof

\brem\label{con.r1}  Let $1<m<k$. Note that one of the matrices in Formula \eqref{con.eq2}
has determinant 1. Thus, one can always find a neighbor $C$ of $E$ in $\Gamma$ and
natural $n$ and $l$ such that $kl-mn=1$ and
$\Gamma$ is obtained by inner blowing ups from the dual graph $\Gamma_0$ of the minimal
resolution of the indeterminacy points of
the function $y^l/x^n$ with the proper transform of  $C$  being the only
$(-1)$-vertex in $\Gamma_0$. 
Actually, the same claim remains true in the case of $m=1$ and $k\geq 1$ if one lets $n=k-1$ and $l=1$
where in the case of $k=m=1$ one has to use the proper transform of the $y$-axis instead of $C$.
\erem 

\blem\label{con.l2} Let the notation of Remark \ref{con.r1} hold. Then 
the functions
$x^k/y^m$ and $y^l/x^n$
form a coordinate system in a neighborhood of the point $C\cap E$
with $C$ and $E$ being the coordinate axes.
\elem

\bproof As we stated in  Lemma \ref{con.l1} either
(a) the functions $x^k/y^m$ and $y^l/x^n$ or (b)  the functions
$y^m/x^k$ and $x^n/y^l$ form a coordinate system in a neighborhood
of $C\cap E$ with $C$ and $E$ being the coordinate axes.
 However, $(y^m/x^k)^l(x^n/y^l)^m =x^{-1}$.
Since the lift of $x$ to a neighborhood of $C\cap E$ is regular we see that (b) is impossible which
yields the desired conclusion.
\eproof

\blem\label{con.l3} Let $\Gamma$  be a connected contractible (not necessarily linear) weighted graph
such that it has only one $(-1)$-vertex $E$.
Then $\Gamma$ contains $n-1$ branch points $E_1, \ldots, E_{n-1}$ such that

{\rm (1)}  each of them has exactly three neighbors;

{\rm (2)} every $E_i$ is contained in a linear subgraph $\Gamma_i$ of $\Gamma$
such that after changing the weight of $E_i$ to -1 this subgraph becomes a weighted graph
 from Lemma \ref{con.l1} with $m,k\geq 2$;

{\rm (3)} every $E_i$ has a neighbor which is an endpoint of $\Gamma_{i+1}$ where
$\Gamma_n$ is as in Lemma \ref{con.l1} (in particular, it contains $E$).

\elem

\bproof We use induction by $n$. If $n=1$, then $\Gamma$ is linear and the conclusion follows from Lemma \ref{con.l1}.
Assume that the statement is true for $n-1$ and we want to prove it for $n$. To get $\Gamma$ contracted
one needs to contract $\Gamma_n$ since otherwise the branch point $E_{n-1}$ prevents 
further contraction. By Lemma \ref{con.l1} $\Gamma_n$ is of the desired form. Contracting $\Gamma_n$
one changes the weight of $E_{n-1}$ and this weight must become -1 to proceed with contraction.
Hence, the induction assumption yields the desired conclusion.
\eproof

\bprop\label{con.p1} Let  $\psi: X \rightarrow X_0$ be a  proper birational morphism of smooth complete
rational surfaces and  $H\subset X$ be a smooth irreducible curve such that $q_0=\psi (H)$ is a singleton. 
Then $\psi =\mu \circ \theta$ where $\theta : X \rightarrow \brX$ is a morphism into
a  smooth complete rational surface $\brX$ and
 $\mu : \brX \to X_0$ is a morphism
such that  the restriction of $\mu$ over $\brX\setminus q_0$ is an isomorphism
and the dual graph $\Gamma$ of  $G=\mu^{-1} (q_0)$
is of the same form as in Lemma \ref{con.l3} with  $\brH=\theta (H)$  being the only $(-1)$-vertex.
Furthermore, $\mu$ can be presented as a composition 
of morphisms $\mu_{i+1} : X_{i+1}'\to X_{i}', \, i=0, \ldots, n-1$ with $\brX=X_n'$ and $X_0'=X_0$  such that
each $\mu_i$ has a connected exceptional divisor whose dual graph is $\Gamma_i$.  

\eprop

\bproof 
Note that $X$ is obtained by  a sequence of blowing up of a finite set in $X_0$ and infinitely near points.
Choose in this sequence the smallest subsequence $R$ which generates a proper transform $\brH$ of $H$
in a surface $\brX$ over $X_0$. That is, $\psi$ factors through a proper birational morphism $\mu : \brX\to X_0$ with
a connected contractible exceptional divisor $G$ which contains $\brH$ as an irreducible component.
Now the first statement is the consequence of the following.

{\em Claim.} The dual graph $\Gamma$ of $G$ is a contractible graph with $\brH$ being a unique $(-1)$-vertex.
Indeed, every $(-1)$-vertex in $\Gamma$ appears in this sequence $R$
of blowing ups only after the proper transforms of its neighbors have been already created (otherwise,
the weight of this vertex cannot be $-1$).
This implies that $\Gamma$ cannot
contain linear $(-1)$-vertices except for $\brH$ since otherwise $R$
is not the smallest subsequence. Hence,  $\Gamma$ satisfies the
assumptions of Lemma \ref{con.l3} and we have the first statement.

The second statement is automatic for $n=1$. Hence, in the general case contracting the components
corresponding to the vertices
of $\Gamma_n$ in $\brX$ we get the desired conclusion by induction.
\eproof
\section{Local coordinate systems}\label{loc}
 The following fact is well-known and the author does not know to whom
its original proof has to be attributed.

\blem\label{loc.l0} Let $X_0$ be a smooth complete rational surface and $q_0\in X_0$.
Then there is a neighborhood of $q_0$ in $X$ isomorphic to the affine plane.
\elem

\bproof
Recall that a variety  $Y$ belongs to class $\cA_0$ if every point in $Y$ has a neighborhood isomorphic to an affine
space (see, \cite[Definition 2.3]{For06}). Any variety obtained by blowing up of $Y$ at a point is still in class $\cA_0$ \cite[Section 3.5D]{Gro89}. In particular, every
smooth complete rational surface  $X$ belongs to class $\cA_0$ since such surface
is obtained from $\PP^2$ or a Hirzebruch surface by consequent blowing up of points  
\eproof

\blem\label{loc.l1} Let $\pi : X \to X_0$ be a proper morphism of smooth projective rational surfaces
such that its restriction is an isomorphism
over $X_0\setminus q_0$ for some $q_0 \in X_0$ and the weighted dual graph $\Gamma_1$ of $\pi^{-1} (q_0)$ is 
the  graph of the minimal resolution of indeterminacy point of a rational function of the form $x^k/y^m$
where $m$ and $k$ are relatively prime. 
Let $C_n$ be the vertex in $\Gamma_1$ that is the proper transform of the exceptional divisor $\tC$
of the blowing up $\tau : \tX\to X_0$ of $X_0$ at $q_0$.
Let $U_0\simeq \A_{x_0,y_0}^2$ be an open subset of $X_0$ with the origin at $q_0$ and
$H$ be the proper transform in $X$ of the $x_0$-axis $H_0\subset U_0$.
Suppose that the graph  $\Gamma_1^+$  of  $\pi^{-1}(H_0)$ is of the form 
$$\cou{H}{}\lin\cou{C_1}{}\lin\cou{C_2}{}\lin \ldots \lin\cou{C_n}{},$$
where $C_1, \ldots, C_n$ are the irreducible components of
$\pi^{-1} (q_0)$.
 Then $\pi$ coincides with  the minimal resolution of indeterminacy of
 $x_0^k/y_0^m$ at $q_0$ and $m\leq k$.
\elem

\bproof
The graph of $\tau^{-1} (H_0)$ is a linear graph consisting of $\tC$ and the proper transform $\tH$ of $H_0$.
To obtain $\Gamma_1^+$ one has to use only inner blowing ups in this graph
starting with the blowing up of the point $\tH\cap \tC$. This implies that $\pi$ is uniquely determined
by $\Gamma_1$.
Hence, $\pi$ coincides either with the minimal resolution of
indeterminacy of  $y_0^m/x_0^k$ or of  $y_0^k/x_0^m$ where $m\leq k$. However, in the latter case 
$H$ meets $C_n$, whereas it is supposed to meet $C_1$.
Thus,  $n=1$ which implies $m=k=1$ and concludes the proof.
\eproof

\bdefi\label{loc.d1} (1)  For $c \in \bk$ and $(x_0,y_0)$ as in Lemma \ref{loc.l1}  a local coordinate system
of the form $(x_0+c,y_0)$ will be called
a {\em bottom coordinate system} for $\pi : X\to X_0$. 

(2)  By Lemma \ref{con.l2} for some relatively prime
$n$ and $l$ the lift of  $x_0^{k}/y_0^{m}$
and  $y_0^{l}/x_0^{n}$ to $X$ yields a coordinate system $(x_1,y_1)$ on an open subset $U_1\simeq \A^2$
of $X$  such that $y_1$ vanishes on the only $(-1)$-vertex of $\Gamma_1$.
We call $(x_1+c,y_1)$ a {\em top coordinate system} for $\pi : X\to X_0$. 
\edefi

\brem\label{loc.r1}  (i)
There is some freedom in the choice of a bottom coordinate system. Indeed,
if one replaced $y_0$ by $y_0+ax_0^l$ where $a \in \bk$ and $l$ is sufficiently large, then
the proper transform of the new $x_0$-axis will still meet $C_1$. Similarly, one can replace $x_0$ by 
$x_0+by_0, \, b \in \bk$. Hence, we can assume that every given finite subset $M$ of $U_0\setminus \{ q_0\}$
does not meet the coordinate axis. 

(ii) A neighborhood $U_0\simeq \A^2_{x_0,y_0}$ of $q_0$
as in Lemma \ref{loc.l1} always exists whenever $k$ and $m$ are relatively prime. 
Indeed, a neighborhood $U_0$ isomorphic to an affine plane always exists by Lemma \ref{loc.l0}.
Consider now a smooth curve $H'\subset X$ meeting
$C_1$ transversely. Then the proper transform $H_0' \subset X_0$ of $H$ contains $q_0$ and it is smooth
since blowups 
preserve SNC type of curves.  Choosing a coordinate system $(x_0,y_0)$ on $U_0$
such that $H_0$ is tangent to $H_0'$ at $q_0$ with sufficiently high order we see that the graph
of $\pi^{-1} (H_0)$ has the desired form.
Using the similar argument we can always find a smooth affine line $F_0\subset X_0$ such
that the graph of $\pi^{-1} (F_0)$ is linear with the proper transform of $F_0$ meeting $C_n$.
\erem

\blem\label{loc.l2} 
Let $\pi : X \to X_0$ be a proper morphism of smooth projective rational surfaces
such that its restriction is an isomorphism
over $X_0\setminus q_0$ for some $q_0 \in X_0$ and the dual graph $\Gamma_1$ of $\pi^{-1} (q_0)$ is 
the graph of a minimal resolution of indeterminacy point of a rational function 
of the form $x^k/y^m$ where $k$ and $m$ are relatively prime.   
Let $U_0\simeq \A_{x_0,y_0}^2$ be an open subset of  $X_0$ with the origin at $q_0$
 such that the graph of $\pi^{-1} (H_0)$ is
linear where $H_0$ is the $x_0$-axis.
Then $\pi$ can be presented as a composition 
of morphisms $\nu_{i} : X_{i}\to X_{i-1}, \, i=1, \ldots, N$ with $X=X_N$  such that

{\rm (i)}  the restriction of $\nu_{i}$ is an isomorphism over $X_{i-1}\setminus \{ q_{i-1}\}$ for some $q_{i-1}\in X_{i-1}$
and the exceptional divisor $\nu_{i}^{-1} (q_{i-1})$ contains  only one $(-1)$-curve $H_i$;

{\rm (ii)} there is an open subset  $U_i\simeq \A^2_{x_i,y_i}$ of $X_i$ such that $q_i\in H_i$
and  $y_i$ vanishes on $H_i\cap U_i$ for $i<N$;

{\rm (iii)} $(x_i,y_i)$ is a bottom coordinate system for  $\nu_{i+1}$ and
a top coordinate system for  $\nu_{i}$; i.e., $\nu_{i+1}$ is
the minimal resolution of indeterminacy points of a function of
the form $(x_i-c_i)^{k_i}/y_i^{m_i}$ where $c_i$ is the $x_i$-coordinate of  $q_i$ 
and, furthermore,  $m_i=k_i=1$ for $i\leq N-2$, while $m_{N-1}\leq k_{N-1}$;

{\rm (iv)} $H_N$ coincides with the only  irreducible component $E$ of $\pi^{-1} (q_0)$ with selfintersection $-1$.

\elem

\bproof Let $C_1, \ldots, C_n$ be as in Lemma \ref{loc.l1}. If 
the proper transform $H$  of $H_0$ meets $C_1$, then the conclusion follows from Lemma \ref{loc.l1} with $\pi=\nu_1$.
In particular, this is true when the length $n$ of $\Gamma_1$ is 1, i.e., $C_1=C_n$ and $\pi$ 
is the blowing up of $X_0$ at $q_0$.
Hence, we suppose that $H$ meets $C_n$ and use induction on $n$.  
Let $\nu_1: X_1 \to X_0$ be the blowing up of $q_0$, i.e., $H_1=\nu_1^{-1} (q_0)$ is the proper transform of $C_n$. 
Then $(x_0,y_0)$ is the bottom coordinate system for $\nu_1$
and $\pi =\nu_1\circ \pi_1$ where the restriction of $\pi_1: X\to X_1$ is an isomorphism over $X_1\setminus\{q_1\}$
for some $q_1\in H_1$.  Note that we have a top coordinate system $(x_1,y_1)$ for $H_1$ with origin at $q_1$
(in particular, $y_1$ vanishes on $H_1$). Note also  that the graph $\Gamma'$
obtained by removing $C_n$ and its adjacent edge from $\Gamma_1$ is a contractible graph
which unique $(-1)$-vertex is still $E$. Furthermore, this graph is the dual graph of $\pi_1^{-1} (q_1)$.
Since its length is $n-1$ we are done by the induction assumption.
\eproof

\bprop\label{loc.p1} Let  $\mu : \brX \to X_0$, $q_0\in X_0$, and $\mu_{i+1} : X_{i+1}'\to X_{i}', \, i=1, \ldots, n-1$ 
be as in Proposition \ref{con.p1}. 
Let $(x_0, y_0)$ be a bottom coordinate system on $\A^2\simeq U_0\subset X_0$ for $\mu_1 : X_1' \to X_0$
with the origin at $q_0$.
Then $\mu$ can be presented as a composition 
of morphisms $\nu_i : X_i\to X_{i-1}, \, i=1, \ldots, N$ with $\brX=X_N$ such that

{\rm (i)}  the restriction of $\nu_{i}$ is an isomorphism over $X_{i-1}\setminus \{ q_{i-1}\}$ for some $q_{i-1}\in X_{i-1}$
and the exceptional divisor of  $\nu_i$ contains  only one $(-1)$-curve $H_i$;

{\rm (ii)} there is an open subset  $U_i\simeq \A^2_{x_i,y_i}$ of $X_i$ such that $q_i\in U_i\cap H_i$
and  $y_i$ vanishes on $H_i\cap U_i$ for $i<N$; 

{\rm (iii)} each $\nu_{i+1}$ is
the minimal resolution of indeterminacy points of a function of
the form $(x_i-c_i)^{k_i}/y_i^{m_i}$ where  $c_i$ is the $x_i$-coordinate of $q_i$  and $m_i\leq k_i$.

\eprop

\bproof   If $n=1$, then we are done by Lemma \ref{loc.l2}.
Otherwise,
 let $\nu_1=\mu_1$ and consider the top coordinate system $(x_1,y_1)$ of $\mu_1$ associated with $(x_0,y_0)$.
In particular, $y_1$ vanishes on the $(-1)$-vertex $H_1'$ of the exceptional divisor of $\mu_1$.
Let  $\Gamma$ and $\Gamma_i$ be as in Proposition \ref{con.p1}.
The  smallest subgraph of $\Gamma$ consisting of the proper transform of $H_1'$ and all vertices
of $\Gamma_2$ is linear. Thus,  we can decompose $\mu_2$ as in Lemma \ref{loc.l2} 
 using $(x_1,y_1)$ as the bottom coordinate system in the first factor of the decomposition.
Now the induction yields the desired conclusion.
 \eproof
\blem\label{loc.l3} Let the assumptions of Proposition \ref{loc.p1} hold, $H_0$ be the $x_0$-coordinate in $U_0$,
and $\ddU_0\simeq \A^2_{\ddx_0,\ddy_0}$ be an open subset of $X_0$ with
the origin at $q_0$ and $\ddy_0$ vanishing on $H_0\cap \ddU_0$, i.e., $(\ddx_0,\ddy_0)$
is another bottom coordinate system for $\nu_1: X_1\to X_0$.
Let $ \A^2_{\ddx_i,\ddy_i}\simeq\ddU_i \subset X_i$ and $\ddc_i$ have the same meaning for
$(\ddx_0,\ddy_0)$ as $U_i$ and $c_i$ for $(x_0,y_0)$. Let $\alpha$  be an automorphism of
$X_0$  such that $\alpha (q_0)=q_0$, $\alpha (U_0)=\ddU_0$, and
$\alpha^* ((\ddx_0,\ddy_0)) = (x_0(1+ p(x_0,y_0)), y_0(1+q(x_0,y_0)))$
where $p$ and $q$ are rational functions on $U_0\simeq \A^2_{x_0,y_0}$ that are  regular
in a neighborhood of the origin  and vanish  at the origin.
Then   for every $i\geq 1$ one has $\ddc_i=c_i$ and 
there exists an automorphism $\alpha_i$ of $X_i$ such that $\alpha_i (U_i)=\ddU_i$ and
$\nu_{i}\circ \ldots \circ \nu_1 \circ \alpha_i= \alpha \circ \nu_{i}\circ \ldots \circ \nu_1$.
\elem
\bproof  Recall that $x_1=x_0^k/y_0^m, y_1=y_0^l/x_0^n, \ddx_1=\ddx_0^k/\ddy_0^m$ and 
$ \ddy_1=\ddy_0^l/\ddx_0^n$ where $kl-mn=1$. 
In particular, 
$(x_0,y_0)=(x_1^ly_1^m, x_1^ny_1^k)$. 
Then $$\alpha^* (\ddx_1)=x_1(1+f_1(x_0,y_0))=x_1(1+f_1(x_1^ly_1^m,x_1^ny_1^k))=x_1(1+p_1(x_1,y_1))$$ and 
$$\alpha^* (\ddy_1)=y_1(1+g_1(x_0,y_0))=y_1(1+g_1(x_1^ly_1^m,x_1^ny_1^k))=y_1(1+q_1(x_1,y_1))$$
where  $f_1$ and $g_1$ are rational functions regular at the origin and vanishing at it.
Since $k,m>0$ we see that $p_1$ and $q_1$ are rational functions regular on the $x_1$-axis  in
 $U_1\simeq \A^2_{x_1,y_1}$  and vanishing on it. 
Hence,  $\ddc_1=c_1$ and for the lift $\alpha_1 : X_1\to X_1$ of $\alpha$ one has $\alpha_1( U_1)=\ddU_1$ and
$$\alpha_1^*((x_1,y_1)) = (x_1(1+p_1(x_1,y_1)), y_1(1+q_1(x_1,y_1))).$$
Let $(\tx_1, \ty_1)=(x_1-c_1, y_1)$ and $(\tilde \ddx_1,\tilde \ddy_1)=(\ddx_1-c_1, \ddy_1)$.
 Since both $p_1(x_1,y_1))$ and $q_1(x_1,y_1))$
vanish on the $x_1$-axis we see that
$\alpha^* (\tilde \ddx_1)=\tx_1(1+\tp_1(\tx_1,\ty_1))$ and $\alpha^* (\tilde \ddy_1)=\ty_1(1+\tq_1(\tx_1,\ty_1))$
where $\tp_1$ and $\tq_1$ vanish at the origin. Thus, replacing $(x_0,y_0)$ and $(\ddx_0,\ddy_0)$
with $(\tx_1,\ty_1)$ and $(\tilde \ddx_1,\tilde \ddy_1)$, respectively,  and using the induction on $N$ we
get the desired conclusion.
\eproof
As an immediate consequence of Lemma \ref{loc.l3} we have the following.
\bprop\label{loc.p2} Let the assumptions of Lemma \ref{loc.l3} hold, 
$R$ be the closure of $U_0\setminus \mu (U_N)$ and $\ddR$ be
the closure of $\ddU_0\setminus \mu (\ddU_N)$. Then $\alpha (R) =\ddR$.
\eprop
\section{The case of surfaces}\label{it}

Let us fix notations for the rest of this section.

\bnota\label{it.n1}
(1) Let $Q'=\PP^1\times \PP^1$ be the quadric that is the natural completion of $\A^2_{x,y}$.
Let $E_b'=\{ y=b\}$ and  $F_a'=\{ x=a\}$  where $a,b \in \PP^1$. Let $q_b'=F_b'\cap E_0'$ and $p_b'=F_b'\cap E_\infty'$.
By $Q_*'\subset Q'$ we denote the union of $\A^2_{x,y}$ and the germ $(Q', p_\infty')$.
We consider also the quadric $Q''$ that is the natural completion of $\A_{s,t}^2$
and suppose that $E_b''$,   $F_a'', Q_*''$, $q_b''$ and $p_b''$ have the same meaning for $Q''$
as $E_b'$,   $F_a',Q_*'$, $q_b'$ and $p_b'$ for $Q'$.  

(2) Let $m$ and $k$ be relatively prime natural numbers
with $m\leq k$. We suppose that $n$ and $l$ are as in Remark \ref{con.r1}.
In particular,  $kl-mn=1$. 
By $\chi : Q'' \dashrightarrow Q'$ we denote the rational map that extends the morphism
$\A^2_{s,t} \to \A^2_{x,y}, \, (s,t) \mapsto (s^lt^m, s^nt^k)$ (so  $\chi (s,t) = (st, t)$ for $m=k=1$).
\enota
\blem\label{it.l1} Let  $c\in \bk$  and $\alpha$
be the automorphism of $Q'$ whose restriction to $\A_{x,y}^2$ is given by $(x,y) \mapsto (x+c,y)$.
Then $\alpha \circ \chi|_{Q_*''}$ yields a morphism 
$\chi_c : Q_*''\to Q_*'$.
\elem
\bproof
Since the restriction of $\alpha$ yields an automorphism of $Q_*'$ it suffices to consider the case $c=0$.   
Note that $p_\infty'$ is contained a neighborhood isomorphic to $\A_{\hx,\hy}^2$
where $(\hx,\hy)=(1/x,1/y)$ and the union of this neighborhood with $\A_{x,y}^2$ contains $Q_*'$.
There is the similar neighborhood $ \A_{\hat s,\hat t}^2$ of $p_\infty''$ whose union with $\A_{s,t}^2$ contains $Q_*''$
where $(\hat s,\hat t)=(1/s, 1/t)$. The formula for $\chi$  in Notation \ref{it.n1} implies that $\chi$ is regular
on $ \A_{\hat s,\hat t}^2$ and
$\chi( \A_{\hat s,\hat t}^2)\subset \A^2_{\hx,\hy}$ (more precisely, $\chi|_{ \A_{\hat s,\hat t}^2 }$ is given by
$(\hat s, \hat t)\mapsto (\hat s^l\hat t^m, \hat s^n\hat t^k)$).  Thus, $\chi$ is regular on $Q_*''$ and $\chi(Q_*'')\subset Q_*'$
which concludes the proof.
\eproof
\blem\label{it.l2} Let $\tau : \brQ \to Q_*'$ be the minimal resolution of indeterminacy
points of the function $(x-c)^k/y^m$ at $q_c'$ and
$\brE$ be the only $(-1)$-component in $\tau^{-1} (q_c')$. 

{\rm (1)} Let $k>1$. Then there is  a morphism
$\breve \chi_c : Q_*''\to \brQ$ such that $ \chi_c =\tau \circ \breve \chi_c$ and
the restriction of $\breve \chi_c$ to $\A^2_{s,t}$ is an embedding.
Furthermore, $\breve \chi _c (E_0''\setminus \{q_\infty''\})\subset \brE$ and
 there is an irreducible component $C$ of $\tau^{-1} (q_c')$ meeting $\brE$
such that $\breve \chi_c (F_0''\setminus \{ p_0''\})\subset  C$.

{\rm (2)} If $k=m=1$, then the same conclusion holds if one replaces $C$ with the proper transform
of $F_c'$ in $\breve Q$.
\elem

\bproof  Let $\breve \chi_c : Q_*'' \dashrightarrow \brQ$ be the lift  of $\chi_c$. 
Since the restriction of $\tau$ over  $Q_*' \setminus \{q_c'\}$ is an isomorphism
we have the regularity of $\breve \chi_c$ over $\tau^{-1} (Q_*' \setminus \{q_c'\})$.
Thus, $\breve \chi_c$ is regular in a neighborhood of $p_\infty''$ since $\chi_c (p_\infty'') \in Q_*' \setminus \{q_c'\}$.
Note that $s=\chi_c^*((x-c)^k/y^m)$  and $t=\chi_c^*(y^l/(x-c)^n)$.
Hence, (1) follows from Lemma \ref{con.l2}.
The similar argument works for (2) with $s=\chi^*((x-c)/y)$ and $t=\chi^* (y)$. 
\eproof

\blem\label{it.l3} Let
$X$ be a smooth complete rational surface, $U_0\simeq \A^2_{x_0,y_0}$
be   an open subset of $X$, and $q=(c,0) \in U_0$.
Let $\nu : \hX\to X$ be the minimal resolution of indeterminacy of the function $(x_0-c)^k/y_0^m$ at $q$,
 that is, $(x_0,y_0)$ is a bottom coordinate system for $\nu$ and we have a
top coordinate system $(x_1,y_1)$ on an open subset $U_1\simeq \A^2_{x_1,y_1}$ of $\hX$
as in Definition \ref{loc.d1}.
Suppose that   $\phi : Q_*' \to X$ is a morphism such that its
restriction to $\A_{x,y}^2$  yields 
the isomorphism $\A^2_{x,y} \to \A^2_{x_0,y_0}$ for which $ \phi^*(x_0,y_0) =(x,y)$.
Then there is a morphism $\hat \phi: Q_*'' \rightarrow \hX$ such that its restriction
yields an isomorphism  $\A_{s,t}^2\to U_1$ with  $ \hat \phi^*(x_1,y_1)=(s,t)$
 and $\phi (p_\infty')=\nu\circ \hat \phi (p_\infty'')$.
\elem

\bproof  Let the notations of Lemma   \ref{it.l2} hold.
Then there is a commutative diagram\\

 
\begin{center}
  \begin{tikzpicture}    
  \setlength{\unitlength}{1cm}
  \put(-2,0){$Q_*''$}
    \put(-1.5,0.1){\vector(1,0){1}} 
  \put(-1.2,0.3){$\breve \chi_c$}
   \put(-0.4,0){$\breve Q$}
       \put(0.1,0.1){\vector(1,0){1}} 
  \put(0.4,0.3){$\kappa$}
   \put(1.2,0){$\hX$}
 
       \put(-0.25,-0.15){\vector(0,-1){1}} 
   \put(1.4,-0.15){\vector(0,-1){1}} 
 \put(-0.2,-0.65){$\tau$}
  \put(1.45,-0.65){$\nu$}
    
     \put(-0.4,-1.5){$Q_*'$}
      \put(1.2,-1.5){$X.$}
       \put(0.1,-1.4){\vector(1,0){1}} 
  \put(0.4,-1.2){$\phi$}
  \put(-1.6,-0.2){\vector(1,-1){1.1}} 
  \put(-1.1,-0.6){$\chi_c$}
 

    \end{tikzpicture}
  \end{center}  

\vspace{1.5cm}
Recall  that $s=\chi_c^*((x-c)^k/y^m)$  and $t=\chi_c^*(y^l/(x-c)^n)$.
Hence, $s=(\phi \circ \chi_c)^*((x_0-c)^k/y_0^m)$  and $t=(\phi \circ \chi_c)^*(y_0^l/(x_0-c)^n)$.
By Definition \ref{loc.d1} we can suppose that $(x_1,y_1) =\nu^*((x_0-c)^k/y_0^m, y_0^l/(x_0-c)^n)$.
Let $\hat \phi :=\kappa\circ \breve \chi_c: Q_*'' \rightarrow \hX$.
Then the commutativity of the diagram implies that $ \hat \phi^*(x_1,y_1)=(s,t)$.   Similarly,
the fact that $\chi_c(p_\infty'')=p_\infty'$ implies the equality
 $\phi (p_\infty')=\nu\circ \hat \phi (p_\infty'')$ and the desired conclusion.
 \eproof

\bprop\label{it.p1}  Let the assumptions of  Proposition \ref{con.p1} hold with $X_0=Q'$
and $H$ being an irreducible component in a complete (not necessarily connected)
 SNC curve $D\subset X$.
Suppose that $\nu_i: X_i\to X_{i-1}$ and $U_i\simeq \A^2_{x_i,y_i}, i=1, \ldots , N$ are as in Proposition \ref{loc.p1},
so $X_N=\brX$.
 Assume that  $\psi (D)$ is a finite subset  of $U_0$ and
 $q_0=\psi (H)$ is the origin in $\A_{x,y}^2\subset Q'$.
Suppose also that  $(x,y)$ coincides with $(x_0,y_0)$.
Then replacing $H$, if necessary, by another irreducible component of $D$
one can find a morphism $\phi_N : Q_*'' \to X_N$ such that its restriction yields
an isomorphism $\A_{s,t}^2 \rightarrow U_N\simeq \A_{x_N,y_N}^2$
for which $\phi_N^* (x_N,y_N)=(s,t)$, $\mu\circ\phi_N (p_\infty'')=p_\infty'$,
and $\theta (D)$ is the union of the curve $H_N=\theta (H)$ and a finite subset.

\eprop

\bproof Let $\mu^{-1} (q_0)$   contain the proper transform $\tH$ of an irreducible
component of $D$ different from $H$. Contracting $H_N$ and some other irreducible components
in the weighted dual graph $\Gamma$ of $\mu^{-1} (q_0)$  we obtain another contractible
graph $\tilde \Gamma$ with the proper transform of $\tH$ as the only $(-1)$-vertex.
Note that we have a morphism $\tilde \mu : \tX \to X_0$ such that $\phi$ factors through $\tilde \mu$
and $\tilde \Gamma$ is the weighted dual graph of $\tilde \mu^{-1} (q_0)$.
Since the replacement of $H_N$ and $\mu$ by $\tH$ and $\tilde \mu$ reduces
the number of vertices in dual graphs we can asssume from the beginning 
that $H$ is the only irreducible component of $D$  whose proper transform is contained in $\mu^{-1} (q_0)$.
Let $\phi_0 : Q_*' \to  X_0$ be the embedding such that  $\phi_0^*(x_0,y_0)=(x,y)$
 and $\phi_0 (p_\infty')=p_\infty'$.
By Lemma \ref{it.l3} we have morphism 
$\phi_1 : Q_*'' \to X_1$ such that $\phi_1^*(x_1,y_1)=(s,t)$   and $\nu_1\circ \phi_1 (p_\infty'')=p_\infty'$.
Note that by our assumption unless $N=1$ none of the irreducible components
of $\nu_1^{-1} (q_0)$ is the proper transform of an irreducible component of $D$.
Proceeding by induction we get a morphism $\phi_N : Q''\to X_N$
such that none of the irreducible components
of $\mu^{-1} (q_0)$  but $H_N$ is the proper transform of irreducible components of $D$
and also $\phi_N^*(x_N,y_N)=(s,t)$  and $\mu\circ\phi_N (p_\infty'')=p_\infty'$
which is the desired conclusion.
 \eproof
 \bprop\label{it.p2} Let the assumptions of Proposition \ref{it.p1}  hold and 
 $\mu$ be invertible over a neighborhood of a point
 $v_0\in X$ contained in  $\A_{x_0}^*\times \A_{y_0}^*\subset U_0$. 
 Let $\alpha : X_0\to X_0$ be the automorphism induced by a map
 $$\alpha (x_0,y_0)=\Big(\frac{ax_0}{a-x_0}, \frac{by_0}{b-y_0}\Big)=:(\ddx_0, \ddy_0)$$
 where $a$ and $b$ are general elements of $\bk$ (in particular, $v_0$ has
 another neighborhood $\ddU_0\simeq\A^2_{\ddx_0,\ddy_0}$). Then replacing $U_0$ and $(x_0,y_0)$ by $\ddU_0$ 
 and $(\ddx_0,\ddy_0)$ one has the following

 {\rm (i)} $ \phi_N (p_\infty'')$ is a general point of $X_N$;
 
 {\rm (ii)} $v_N=\mu^{-1} (v_0)$ is contained in $\A_{x_N}^*\times \A_{y_N}^*\subset U_N$.
 
 \eprop
 \bproof 
Note that
 $(\ddx_0,\ddy_0)$ is another bottom coordinate system  for $\nu_1: X_1 \to X_0$.
 Let  $\A^2_{\ddx_i,\ddy_i}\simeq \ddU_i$ have the same meaning as in Lemma \ref{loc.l3},
 $R$ be the closure of $U_0\setminus \mu (U_N)$, and $\ddR$ be
the closure of $\ddU_0\setminus \mu (\ddU_N)$.  
By Proposition \ref{loc.p2} $\alpha (R)= \ddR$.
Since $a$ and $b$ are general and $v_0\in \A_{x_0}^*\times \A_{y_0}^*$
we see that $v_0$ does not belong to $\ddR$ which implies (ii).
Since the replacement of $(x_0,y_0)$ by 
 $(\ddx_0,\ddy_0)$ makes $p_\infty'$ a general point of $X_0$ we have also  (i) which concludes the proof.
 \eproof
 In the next lemmas we shall justify the assumptions of Propositions \ref{it.p1} and \ref{it.p2}.

 \blem\label{it.l4} Let $\psi : X\to \PP^n$ be  a birational morphism
of smooth complete varieties, $v \in X$, and $\psi$ be invertible over a neighborhood of
$v_0=\psi (v)$.  Let $\cD$ be a line bundle on $X$ and $D$ be a divisor in $X$
such that $\cD=\cO_X(D)$.
Suppose that  the pushforward of $D$ generates a nontrivial line bundle on $\PP^n$.
Then for an appropriate choice of $D$
there is a birational morphism $\phi: B\times \PP^1\to X$ such that

{\rm (i)}  $\phi$ is invertible over  a neighborhood of $v$;

{\rm (ii)} there is a general point $\fo$ of some irreducible component
of $D$ such that each curve $\phi ( b\times \PP^1), \, b \in B$ meets $D$ at $\fo$ only
and the intersection of this curve and $D$ is transversal;

{\rm (iii)} the restriction of $\cD$ to $ \phi (b \times \PP^1)$ is a nontrivial line bundle.

\elem

\bproof 
Let $L$ be a hyperplane in
$\PP^n$ that does not contain $v_0$. One can choose a rational function
$f$ on $\PP^n$ with a divisor $T$ such that $T + \psi_* D = mL$.
Note that $\cO_X(D)=\cO_X (\psi^*(T)+D)$. Hence, we  can suppose that $\psi_* D =mL$.
By the assumption $m\ne 0$.
Choose a general point $\fo'$ in $L$ and a general hyperplane  $\PP^{n-1}$ in $\PP^n$.
Then for $\A^{n-1}\subset \PP^{n-1}$ one can consider the set of lines each of which passes through $\fo'$
and a point of $\A^{n-1}$. This set generates a morphism $\phi_0 : \A^{n-1}\times \PP^1\to \PP^n$.
Since $\fo'$ is general we can suppose that $\psi$ is invertible over a neighborhood of the line
$\phi_0 (b_0\times \PP^1)$ containing $v_0$. Thus, $\psi$ is invertible over $\phi_0 (B\times \PP^1)$
for some open subset $B\subset \A^{n-1}$ containing $b_0$. 
The lift of $\phi_0|_{B\times \PP^1}$ to $X$ yields $\phi : B\times \PP^1 \to X$ 
such that every curve $\phi (b \times \PP^1)$ meets  $D$ only at a general point $\fo$ in
the proper transform $H$ of $L$. Hence, the restriction of $\cD$ to $\phi (b \times \PP^1)$ 
is a nontrivial line bundle which concludes the proof.
\eproof
\blem\label{it.l5} Let  $P$ be either $\PP^2$ or a Hirzebruch surface $\F_n$
and let $\beta : Y \to P$ be a birational morphism from
a smooth complete surface $Y$ such that it is invertible over a neighborhood of   $u_0\in P$.
Let $M$ be a finite subset of $Y$ and
 $T$ be a complete SNC curve in $Y$ such that $u=\beta^{-1} (u_0) \notin T\cup M$ and $\beta (T)$ is finite.
Then there is a commutative diagram

\begin{equation}\label{it.eq6}
\begin{array}{ccl} X &  \stackrel{{\lambda}}{\rightarrow} & Y\\
\, \, \,  \downarrow^{\psi}   & 
& \,  \downarrow^{\beta}\\
Q' &  \stackrel{{\gamma}}{\dashrightarrow} &   P
\end{array} 
\end{equation}
where $X$ is a smooth complete rational  surface, the morphisms $\psi$ and
$\lambda$ are birational and $\lambda$ is  invertible over a neighborhood
of $\{ u\}\cup T\cup M$.  Furthermore, one can suppose that $\psi : X\to X_0=Q'$
satisfies the conclusions of Proposition \ref{it.p2} with $D=\lambda^{-1} (T)$,  and $v_0=\psi (\lambda^{-1} (u))$.
\elem

\bproof   Let $M_0$ be a finite subset of $P$ such that $u_0\notin M_0, \beta (M)\subset M_0$
and $\beta$ is invertible over $P\setminus M_0$  (note that $\beta (T) \subset M_0$).
For $P=\PP^2$ let $L$ be a general line in $\PP^2$. Then $L$ does not meet $\{u_0\} \cup M_0$.
Thus, $U=\PP^2\setminus L$ is a neighborhood of $u_0$ isomorphic to $\A^2_{x_0,y_0}$.
Let $F$ be any irreducible component of $T$. Then $\beta$ admits a decomposition as in Proposition \ref{con.p1}
with $\psi : X\to X_0$ and $H$ replaced by $\beta : Y \to \PP^2$ and $F$. Consequently,
 $\beta$ admits a decomposition as in Proposition \ref{loc.p1} with $(x_0,y_0)$
having the same meaning as the bottom coordinate system in  Proposition \ref{loc.p1}.
By Remark \ref{loc.r1} the coordinate system $(x_0,y_0)$ can be chosen so that 
 $u_0 \in \A^*_{x_0}\times \A^*_{y_0}$. 
 Choosing $Q'$ as the natural completion of 
$\A_{x_0,y_0}^2$ we  get the birational map $\gamma : Q'\dashrightarrow \PP^2$
which gives rise to commutative diagram in \eqref{it.eq6} with $\lambda$ being a rational map. 
However, resolving the indeterminacy points of $\lambda$ (which are located in  $\beta^{-1} (L)$)
we can suppose that $\lambda$ is a morphism. Let $H$ be the proper
transform of $F$ in $X$.
By construction all the assumptions
of Proposition \ref{it.p1} for $\psi : X\to Q'$, $H, D$, $v_0$, and $(x_0,y_0)$ are valid which yields
in turn the conclusions of Proposition \ref{it.p2}. 

In the case $P=\F_n$ let $L$ be a fiber of the natural projection $\F_n \to \PP^1$
 that is disjoint from $\{ u_0\} \cup M_0$.
 Choose in $\F_n\setminus L$ a neighborhood of $u_0$ isomorphic to $\A^2_{x_0,y_0}$
 such that  $u_0 \in \A^*_{x_0}\times \A^*_{y_0}$. 
 Since there is a rational map $\gamma : Q'\dashrightarrow \F_n$ whose restriction over $\F_n\setminus L$  is an isomorphism, repeating the argument from the previous case we get the desired conclusion.
\eproof
\bprop\label{it.p3} Let the assumptions and conclusions of Lemma \ref{it.l5} hold.
Then there is an irreducible component $H$ of $T$ and a  morphism $\phi : B\times \PP^1 \to Y$
such that

{\rm (i)}  $\phi$ is locally invertible over a neighborhood of $u$;

{\rm (ii)} there is a general point $\fo \in H$  such that for every $b \in B$ the curve $\phi (b\times \PP^1)$ meets 
$T$ at $\fo$ only and the intersection of  $\phi (b\times \PP^1)$ and $H$ is transversal;

{\rm (iii)} the image of $\phi$ does not meet $M$.

\eprop
\bproof  Since $\lambda$ is locally invertible over a neighborhood of $\{ u\}\cup T\cup M$
we can replace $Y, T$ and $u$ with $X, D$ and $v=\lambda^{-1} (u)$ respectively.
Recall that we have a decomposition $\psi =\mu\circ \theta$ as in Propositions \ref{con.p1} and \ref{loc.p1}
where $\theta : X \to \brX=X_N$ and $\mu: X_N \to Q'$. 
Furthermore, there is an irreducible component $H$ of $D$ such that
$\theta (D)$ is contained in  the union of the curve $H_N=\theta (H)$ and the set $M_N$ of points over which
$\theta$ is not invertible. 
By Proposition \ref{it.p1}
we have a morphism $\phi : Q_*'' \to X_N$ which restricts to an isomorphism $\A_{s,t} \to U_N$
onto some open subset $U_N \subset X_N$.
By Proposition \ref{it.p2} we can suppose that $U_N$ contains
$v_N=\mu^{-1} (v_0)=\theta (v)$. 
Since $\psi$ is invertible over $v_0$, $\theta$ is invertible over $v_N$ and,
thus, $v_N$ does not belong to the finite set $M'=M_N\cup \theta (M)$.
Let  $M''=\phi_N^{-1} (M')\cap \A^2_{s,t}$ and $v_N=\phi_N(s_0,t_0)$ where 
by Proposition \ref{it.p2} $(s_0,t_0)\in \A^*_{s}\times \A^*_t$.
Consider all $a,b$ and $c\in \bk$ for which the image of the morphism
 $\gamma_{a,b,c}: \A_z^1 \to Q_*'', \, z \mapsto (a+bz+cz^2,z)$ 
contains $(s_0,t_0)$. 
Note that  for any point $(s_1,t_1)\ne (s_0,t_0)$ in $\A_{s,t}^2$  and $a \in \bk$ one can choose $b$ and $c$ such 
that the image of  $\gamma_{a,b,c}$ does not contain  $(s_1,t_1)$. In particular, we can suppose
that for a general $a_0$ and some pair $(b_0,c_0)$ the image of $\gamma_{a_0,b_0,c_0}$ does not meet $M''$.  
Hence, for an appropriate  neighborhood $B\subset \A^1$ of $b_0$ we have a morphism
$B\times \A^1_z\to \A_{s,t}^2,\, (b,z)\mapsto (a_0+bz+c_0z^2, z)$ invertible over $(s_0,t_0)$ 
such that its image does not meet $M''$.
This map extends to a morphism $\gamma : B\times \PP^1\to Q_*''$
such that $\gamma (B\times \infty)=p_\infty''$. Note that
$\phi_N \circ \gamma (B\times \PP^1)$ does not meet $M'$
since the restriction of $\phi_N$ to $\A_{s,t}^2$ is an isomorphism onto $U_N$ and
$\phi_N(p_\infty'')$ is a general point of $X_N$ by Proposition \ref{it.p2}.
Furthermore, every curve $\phi_N\circ \gamma (b\times \PP^1)$ meets $H_N$ transversely
and at one general point only since every curve $ \gamma (b\times \PP^1)$ meets
the $s$-axis transversely and only at the point with $s$-coordinate $a_0$.  Since the image
of $\phi_N \circ \gamma$  does not meet $M_N$ we see that
$\phi_N \circ \gamma$ admits a lift $\phi : B\times \PP^1 \to X$ such that (i)-(iii) hold. Hence, we are done.
\eproof
\blem\label{it.l6} Let $\pi : \tX\to X$ and $\tilde \phi : B\times \PP^1 \to \tX$ be
 birational morphisms of smooth varieties such $\pi$ is locally invertible over a neighborhood of $v \in X$ 
 and $\tilde \phi$ is locally invertible
 over a neighborhood of $\tv=\pi^{-1}(v)$.
Let  $\cD$ be a line bundle on $X$ and $\tilde \cD=\pi^* \cD$.
Suppose that  for $b\in B$ such that $\tv \in \tilde \phi (b \times \PP^1)$ the restriction of $\tilde \cD$
to $\tilde \phi (b \times \PP^1)$ is a nontrivial line bundle. Then there exists 
$\phi : B\times \PP^1 \to X$ locally invertible over a neighborhood of $v$ such that 
the restriction of $\cD$
to $\phi (b \times \PP^1)$ is a nontrivial line bundle. 
\elem

\bproof Put $\phi=\pi\circ \tilde \phi$.
Note that if the restriction of $\cD$ to $\phi(b \times \PP^1)$ is
a trivial line bundle, then so is the restriction of $\tilde \cD$ to $\tilde \phi(b \times \PP^1)$. 
This yields the desired conclusion.
\eproof
\bthm\label{it.t1} 
 Let $X$ be a smooth complete rational surface and  $\cD$ be a nontrivial line bundle on $X$
with zero section $Z$.
Then $\cD\setminus Z$ is elliptic. 
\ethm

\bproof By Lemma \ref{loc.l0} $X$ is of class $\cA_0$. 
Hence, $X$ is elliptic since every variety of class $\cA_0$ is elliptic \cite[Section 3.5]{Gro89}.
Thus, by Corollary \ref{pre.c1}  it suffices to show that 
for every $v\in X$ one can find a birational morphism $\phi : B\times \PP^1 \to X$ invertible over
a neighborhood of $v$ such that   for $b_0 \in B$ with $v \in \phi (b_0 \times \PP^1)$
the restriction of $\cD$ to $\phi (b_0 \times \PP^1)$ is a nontrivial line bundle.
Let  $U\simeq \A^2$ be a neighborhood of $v$ in $X$. Letting $\PP^2$ be a completion
of $U$ we have a rational map $\beta: X\dashrightarrow \PP^2$ which is regular and invertible over $U$.
Resolving  indeterminacy points of this map we get a morphism $ \tX\rightarrow \PP^2$ invertible over $U$.
By Lemma \ref{it.l6} we can suppose that $X=\tX$, that is $\beta$ is regular. By Lemma \ref{it.l4} we can assume
that $\beta (\supp D)$ is finite for some SNC divisor $D$ such that $\cD =\cO_X(D)$. 
By Proposition \ref{it.p3} there is
a birational morphism $\phi : B\times \PP^1 \to X$
such that the curve $\phi (b_0\times \PP^1)$ passing through $v$ meets  only one
irreducible component $H$ of $D$. Hence, the restriction of $\cD$
to $\phi (b_0\times \PP^1)$ is a nontrivial line bundle which yields the desired conclusion.
\eproof
\section{Main Theorem}\label{main}

The aim of this section is the following.

\bthm\label{main.t1}
Let $X$ be a complete  uniformly 
rational variety and $\cD\to X$ be a nontrivial  line bundle on $X$
with zero section $Z_\cD$. Then $Y=\cD\setminus Z_\cD$ is elliptic.
\ethm

By Corollary \ref{pre.c1} it suffices to show that 
for every $v\in X$ one can find a smooth variety $B$ and a birational morphism $\phi : B\times \PP^1 \to X$ invertible over
a neighborhood of $v$ such that   for $b_0 \in B$ with $v \in \phi (b_0 \times \PP^1)$
the restriction of $\cD$ to $\phi (b_0 \times \PP^1)$ is a nontrivial line bundle.
To establish this fact we need to check the following property for every $v\in X$.

\bdefi\label{main.d1} Let $X$ be a smooth  complete rational variety, 
$v$ be a point in $X$ with
a neighborhood isomorphic to an open subset of $\A^n$ (which is automatic in the
case when $X$ is uniformly rational),
and $\cD$ be a  nontrivial line bundle on $X$. 
We say that the triple $(X, \cD, v)$ has $\cP$-property
if  
either

\noindent (a) \hspace{0.3cm} there is  an SNC divisor $D'$ such that $\cD=\cO_X (D')$  and
 a morphism of $X$ onto $\PP^n$ such that
its restriction to a neighborhood $W$ of $v$ is an embedding and the pushforward of $D'$ is not a principal divisor on $\PP^n$, or,

\noindent (b) \hspace{0.3cm} 
  given an SNC divisor $D'$ such that $\cD=\cO_X (D')$
 and a morphism of $f: X\to\PP^n$ such that
its restriction to $W$ is an embedding and $\dim f(D') \leq \dim X-2$, 
there is a birational morphism $\phi : B\times \PP^1 \to X$ invertible over
a neighborhood of $v$ such that  $\phi (B \times \PP^1)$ meets $D'$
at a general point $\fo$ of one irreducible component $H$ of $D'$ only
 and the intersection of every curve $\phi (b \times \PP^1), \, b\in B$ and
$H$ at this point is transversal.
\edefi

In fact, we do not need to establish this property for every triple as in Definition \ref{main.d1} since by
Lemma \ref{it.l6}  we can replace $(X,\cD,v)$ by $(\tX, \tilde \cD, \tv)$
where $\sigma : \tX \to X$ is a proper birational morphism from a smooth variety $\tX$
invertible over a neighborhood of $v$, $\tilde \cD=\sigma^*\cD$, and $\tv =\sigma^{-1} (v)$.
Furthermore, by  Lemma \ref{it.l4} and Corollary \ref{pre.c1}, Theorem \ref{main.t1} is valid under 
the assumption (a) of Definition \ref{main.d1}. Thus, we suppose further that
 there is a  birational morphism $\theta : X\to \PP^n$
invertible over a neighborhood of $\theta (v)$ such that $\theta(\supp D')$ is contained in
the indeterminacy set of $\theta^{-1}$. 

\blem\label{main.l1} 
There exists a proper birational morphism  $\sigma : \tX \to X$  as above
such that  replacing  $(X,\cD,v)$ by $(\tX, \tilde \cD, \tv)$ one has 
 a $\PP^1$-bundle $\pr : Q\to S$
over $S\simeq \PP^{n-1}$ and a proper birational morphism
$\alpha: X\to Q$ for which

{\rm (i)} $\alpha (\supp D')$ is contained in the
indeterminacy set  $P$ of $\alpha^{-1}$
and none of the fibers of $\pr$ is contained in $P$;

{\rm (ii)} there is a neighborhood $W\subset S\setminus \pr (P)$ of $s_0=\rho (v)$
for which 
$\rho|_{\rho^{-1} (W)} : \rho^{-1} (W) \to W$ is a $\PP^1$-bundle
where $\rho =\pr \circ \alpha : X \to S$;

{\rm (iii)} $\pr$ has two disjoint sections $\cS$ and $\cS_0$.
\elem

\bproof  Let $P$ be the indeterminacy set of $\theta^{-1}$.
We consider a general line $L$  through $v_0=\theta (v)$ in $\PP^n$ (in particular, $L$ 
does not meet $P$),
a general point $p_0 \in L$, and a general hyperplane $S\simeq \PP^{n-1}$  through $v_0$ in $\PP^n$.
Let $\gamma: Q\to\PP^n$ be  the blowing up of $\PP^n$ at $p_0$.
Then one has the natural morphism $\pr: Q\to S$ which is a $\PP^1$-bundle.
By Lemma \ref{it.l6}  we can replace $X$ with  its blowing up $\brX$ at $\theta^{-1} (p_0)$
(and we also replace $v$ by its preimage  in this blowing up). This yields
the morphism $\alpha : \brX\to Q$. 
Note that the indeterminacy set  $\brP$ of $\alpha^{-1}$ coincides with  $\gamma^{-1}(P)$.
Since $L$  does not meet $P$ there is a neighborhood $W\subset S\setminus \pr (\brP)$ of $\pr(v_0)$ such that
 $\pr^{-1} (W)\cap\brP=\emptyset$.
 Consequently,  $\rho|_{\rho^{-1} (W)} : \rho^{-1} (W) \to W$ is a $\PP^1$-bundle.
  Since $p_0\notin  P$
 the exceptional divisor of $\gamma$ does not meet  $\brP$.
 Hence, none of the fibers of $\pr$ is contained in $\brP$ because
 the exceptional divisor is a section  $\cS_0$ of the $\PP^1$-bundle 
 (the other section $\cS$ is generated by $S$). 
This concludes the proof.
\eproof

\blem\label{main.l2}
There exists a proper birational morphism  $\sigma : \tX \to X$  as before Lemma \ref{main.l1}
such that  replacing  $(X,\cD,v)$ by $(\tX, \tilde \cD, \tv)$ one has 
a $\PP^1$-bundle
$ \pr : Q \to S$ over a smooth complete variety $S$ and
a proper birational morphism $ \alpha : X\to  Q$ for which

{\rm (i)} $\alpha (\supp D')$ is contained in the
indeterminacy set $P$ of $\alpha^{-1}$, $\dim  \rho (H)=n-2$  for  $ \rho= \pr \circ  \alpha : X \to S$ and every
 irreducible component $H$  of $D=\alpha^{-1} (P)$,
and none of the fibers of $\pr$ is contained in $P$;

{\rm (ii)}  there is a morphism $\pi: S\to \PP^{n-1}$ whose restriction to  
a neighborhood  $W\subset S$ of $s _0= \rho (v)$ is an embedding,
the restriction of $ \rho$ over $W$  is $\PP^1$-bundle, and
$ \rho^{-1} (W)\cap D=\emptyset$;

 {\rm (iii)} $\pr$ has two disjoint sections $\cS$ and $\cS_0$.
\elem

\bproof  Let $\alpha: X\to Q$,
 $\rho : X\to S$,  $\pr : Q\to S$, and $P$ be as in Lemma \ref{main.l1}.
Let $\pi : \tS\to S$ be a proper birational morphism
 from a smooth variety $\tS$,
$\tX=X\times_{S}\tS$, and $\tQ=Q\times_{S}\tS$.
Let $\tilde \rho : \tX\to \tS$, $\pr_X : \tX \to X$, $\pr_Q : \tQ \to Q$ be the natural
projections and $\tilde \alpha: \tX\to \tQ$ be the birational morphism induced by $\alpha$.
By Raynaud's flattening theorem  $\pi : \tS\to S$ can be chosen
 such that  $\tilde\rho$ is flat. Furthermore, since $\rho$ is flat over $s_0=\rho (v)$ we can suppose that
$\pi$ is invertible over a neighborhood $W$ of $s_0$ in $S$ \cite[Chapter 4, Theorem 1]{Ray72}.
In particular, the restriction of $\pr_X$ over $\rho^{-1}(W)$ is invertible and Lemma \ref{main.l1} (ii) yields (ii).
Let $\tP$ be the indeterminacy set of $\tilde \alpha^{-1}$, so $\tP\subset \pr_Q^{-1} (P)$.
Since $\tQ$ is smooth every irreducible component $\tH$ of $\tD=\tilde \alpha^{-1} (\tP)$ is a hypersurface in $\tX$
by van der Waerden's theorem (e.g., \cite[Chap. II, Section 7.3]{DaSh94}).
Hence, $\tilde \rho (\tH)$ must have dimension $n-2$
since each fiber of $\tilde \rho$ is one-dimensional because of flatness. 
 By Lemma \ref{main.l1} (i) $\pr^{-1} (s) \setminus P$ is a
nonempty open subset of $\pr^{-1} (s)\simeq \PP^1$ for every $s \in S$. Hence, 
 $\tilde \pr^{-1} (s) \setminus \tP$ is a
nonempty open subset of $\tilde \pr^{-1} (s)$ for every $\ts \in \tS$. 
Thus,  none of the fibers of $\tilde\pr$ is contained in $\tP$.
Thus, we have (i)
which yields the desired conclusion,  since (iii) follows from
Lemma \ref{main.l1} (iii).
\eproof

\blem\label{main.l3} Let the conclusions of Lemma \ref{main.l2} hold and $F$ be an irreducible
component of $T= \rho (D)$.
Suppose that property $\cP$ is true for the triple
$(S,\cO_S (F), s_0)$.  Then one can find a morphism $\phi : B \times \PP^1 \to S$
invertible over a neighborhood of $s_0$
such that 
there is a general point $\fo \in F$ which is also a smooth point
of  every curve $\phi(b\times \PP^1), \, b \in B$ and $F$ meets $\phi(b\times \PP^1)$ transversely at $\fo$.
\elem

\bproof 
If  $F'=\pi(F)$  is a hypersurface in $\PP^{n-1}$, then consider  
a general point $\fo'$ in $F'$, 
a neighborhood $B$
of $\pi(s_0)$ in a  hyperplane $L\simeq \PP^{n-2}$ of $\PP^{n-1}$ containing $\pi(s_0)$ 
but not $\fo'$, and the set of lines passing
through $\fo'$ and the points of $B$.  One can choose $B$ such that this set generates 
a morphism $\psi: B \times \PP^1 \to \PP^{n-1}$.
Since $\fo'$ is a general point  of $F'$  the morphism $\pi$ is invertible
over a neighborhood of the line through $\pi(s_0)$ and $\fo'$.  Hence, for an appropriate $B$
we can suppose that $\psi$ admits a lift to $S$ which is a
desired morphism $\phi : B \times \PP^1 \to S$. If $\dim F'\leq n-3$, then the desired conclusion
follows from the assumption (b) of property $\cP$ and we are done.
\eproof
%
\blem\label{main.l5} 
Let the notation of Lemmas \ref{main.l2} and \ref{main.l3} hold,
 $\hQ=(B\times \PP^1)\times_{S,\phi,\pr} Q$,  $\hat \pr : \hQ \to B\times \PP^1$ be the natural projection,
 and $\hQ^b=\hat \pr^{-1} (b\times \PP^1)$.
Replacing $B$ by an affine neighborhood of a given point $b_0\in B$  
one can find an isomorphism  $\mu : \hQ\to B\times \hQ^{b_0}$
over $B\times \PP^1$.
\elem
\bproof 
By Lemma \ref{main.l2} (iii) $\hat  \pr$ has two disjoint sections $\hat \cS$ and $\hat \cS_0$.
Thus, $\hQ\setminus \hat \cS_0$ can be viewed as a line bundle $\cL$ over $B\times \PP^1$
and $\hQ$ as the projectivization of $\cL\oplus \cE$ where $\cE$ is the trivial line bundle.
Let $\cL_b$ and $\cE_b$ be the restrictions of these bundles to $b\times \PP^1$.
Replacing $B$ by an appropriate affine neighborhood of $b_0$ we can suppose that
a divisor generating $\cL$  is linearly equivalent to a divisor of the form
$m( B\times q)$ where $q\in \PP^1$ and $m\in \Z$. Thus, $\cL=\kappa ^*\cL_{b_0}$
where $\kappa  : B\times \PP^1 \to \PP^1 \simeq b_0\times \PP^1$ is the natural projection.
Consequently,  $\kappa^*(\cL_{b_0}\oplus \cE_{b_0}) =\cL \oplus \cE$ and
$\hQ$ is  isomorphic to $B\times \hQ^{b_0}$ over $B\times \PP^1$.
\eproof
\bnota\label{main.n1} 
Let the notation of Lemmas \ref{main.l3} and  \ref{main.l5} 
hold,  $F^*$ be a smooth affine subset of
$F$ containing $\fo$, and  $V$ be an affine neighborhood of $F^*$ in $S$.
Shrinking $B$ we can   suppose
that $\phi^{-1} (\fo)$ coincides with $R=B\times 0$ where $0\in \PP^1$.
We  also let $\hY_0^{b_0}=\hat \pr^{-1}(\phi^{-1} (V \setminus T) \cup R)\cap \hQ^{b_0}$
and $U=(\phi^{-1} (V \setminus T)\cup R) \cap \hat \pr\circ \mu^{-1} (B\times \hY^{b_0})$
(observe that $\phi^{-1} (F^*)\cap U=R$).

\enota
Restricting the isomorphism $\mu$ from Lemma \ref{main.l5} we get the following.
\blem\label{main.l6}
For $\hY_0=\hat \pr^{-1}(U)$ there is an embedding 
$\nu_0: \hY_0 \hookrightarrow B\times \hY_0^{b_0}$  over $B\times \PP^1$.

\elem
\blem\label{main.la} Let $Y'\overset{\alpha'}{\to} Y_0'\overset{\pr}{\to}Q'$ be projective morphisms of smooth varieties,
 $\dim Y'=\dim Q'+1$, and $\rho=\pr\circ \alpha'$.
 Let $F_0$ be a smooth hypersurface in $Y_0'$ such that $F=\pr (F_0)$
 is a hypersurface in $Q'$ and the restriction of $\alpha'$ over $Y_0'\setminus F_0$ is invertible.
Then one can find an open subset $V$ of $Q'$ with a nonempty smooth $F^*=V\cap F$ 
such that for  $Y=\rho^{-1} (V), Y_0=\pr^{-1} (V)$, and $\alpha=\alpha'|_Y$
one has the following

{\rm (i)} $\alpha : Y \to Y_0$ can be decomposed into projective
birational morphisms $\alpha_i : Y_{i}\to Y_{i-1}, \, i=1, \ldots, m$ of smooth varieties
with irreducible exceptional divisors  $E_i$;

{\rm (ii)} for $C_{i-1}=\alpha_i(E_i)$ the natural projection $C_{i-1}\to F^*$ is an unramified covering;

{\rm (iii)} $\alpha_i$ is the blowing-up of the defining ideal of $C_{i-1}\subset Y_{i-1}$ with
 $E_i\simeq C_{i-1}\times \PP^1$;

{\rm (iv)} the exceptional divisor $E$ of $\alpha$ is of SNC type.

Furthermore, let $G\subset V$ be a curve meeting $F^*$ transversely at a point $\fo$ only.
Suppose that the surface $S_0=\pr^{-1}(G)$ is smooth over $\fo$ and it
meets $C_0$ transversely. 
Then the surface $\rho^{-1}(G)$  is smooth over $\fo$ and it
meets $E$ transversely. 
\elem
\bproof
Let $F^*$ be a nonempty smooth affine open subset of $F$ such that it admits a finite morphism onto
an open subset $W$ of  $\A^{n-2}$. 
Choose an affine neighborhood $V\subset Q'$ of $F^*$ such that
the morphism $F^* \to W$ extends to
a morphism $V\to W$  and, thus, $Y$ and $Y_0$ can be viewed as varieties over $W$. 
 Let  $\bK$ be the algebraic closure of the field of rational functions on $W$.
 For every algebraic variety $X$  over $W$ we denote by $X^\bK$ the variety $\Spec\bK\times_{\Spec \bk} X $ 
and we also use upper index $\bK$ for points in  $X^\bK$ and  functions on it.
 Then we have  the projective morphisms
 $\alpha^\bK: Y_N^\bK\to Y_0^\bK$  
 and $\pr^\bK : Y_0^\bK\to V^\bK$
  where $  Y_N^\bK$ and $Y_0^\bK$ are smooth surfaces, $V^\bK$ is smooth curve, 
 and the restriction of $\alpha^\bK$ is an isomorphism over $Y_0^\bK\setminus F_0^\bK$.
 By \cite[Chap. V, Proposition 5.3]{Har04} $\alpha^\bK$ can be presented as a composition
of monoidal transformations   $\alpha_i^\bK : Y_i^\bK\to Y_{i-1}^\bK, \, i=1, \ldots$ with centers at  
points $q_{i-1}^\bK$  above  $F_0^\bK$.  
Then  $q_0^\bK$ is a generic point of some closed hypersurface  $C'$ in $F_0$ such that $\pr (C')=F$.
Since $\dim F=\dim C'$ the preimage of a general point in $F$ under the morphism $\pr|_{C'}$
is finite. This implies that $F^*$ can be chosen so that
for some smooth open subset $C_0$ of $C'$ the morphism $\pr|_{C_0}: C_0\to F^*$ is finite and, furthermore,
by  \cite[Chapter III, Corollary 10.7]{Har04} we can suppose that $\pr|_{C_0}$  is smooth
and, therefore, $C_0\to F^*$ is an unramified covering. Now $\alpha_1^\bK: Y_1^\bK \to Y_0^\bK$
corresponds to the blowing up  $\alpha_1: Y_1 \to Y_0$ of $Y_0$ 
along the smooth center $C_0$. Hence,  $\alpha_1|_{E_1} : E_1\to C_0$ is a
$\PP^1$-bundle and shrinking $F^*$ again we can suppose that $E_1$ is naturally isomorphic
to $C_0\times \PP^1$.
Using the same argument and induction we see that for an appropriate $F^*$ 
one gets morphisms $\alpha_i: Y_i \to Y_{i-1}$ each of which is the blowing up of $Y_{i-1}$ along smooth $C_{i-1}$
with $E_i \simeq C_{i-1}\times \PP^1$  and $\pr\circ \alpha_1 \circ \ldots \circ \alpha_{i-1}|_{C_{i-1}} : C_{i-1} \to F^*$ 
being an unramified covering.  Thus, we have (i)-(iii).
Statement (iv) follows again from a choice of an appropriate $F^*$ and
the fact of the exceptional divisor of $\alpha^\bK$ is an SNC curve.

Because of transversality the restriction of the defining ideal of $C_0$ in $Y_0$
to $S_0$ is the defining ideal of $\pr^{-1} (\fo)\cap C_0$ in $S_0$. 
Hence,  by \cite[Chapter II, Corollary 7.15]{Har04} 
the surface $S_1=\alpha_1^{-1} (S_0)$ is obtained from $S_0$ by the monoidal
transformations at the points of  $\pr^{-1} (\fo)\cap C_0$. Consequently, $S_1$ is smooth over $\fo$.
Let $C_1'$ be any hypersurface in $E_1$ \'etale over $F^*$  (e.g., $C_1'=C_1$). 
For every point $\fo_1 \in S_1 \cap C_1'$ the map 
$(\pr\circ \alpha_1)_*: T_{\fo_1}C_1' \to T_\fo F^*$ is an isomorphism. 
Since $(\pr\circ \alpha_1)_*(T_{\fo_1}S_1)\subset T_\fo G$ and $T_\fo G \cap  T_\fo F^*=0$  we
have $T_{\fo_1}C_1' \cap T_{\fo_1}S_1=0$ which yields the transversality of the intersection
of $S_1$ and $C_1'$. Therefore, $S_1$ is transversal to $E_1$.
Now the induction implies the second statement and we are done.
\eproof
\bnota\label{main.n2}  We keep Notation \ref{main.n1} and  the notations of Lemma \ref{main.l6}.
 For $Y= \rho^{-1} (V)$ and $Y_0=\pr^{-1} (V)$
 we have the morphisms $\alpha_V : Y\to Y_0$, $\pr_V : Y_0 \to V$, and  $\rho_V : Y\to V$
that are the restrictions of $\alpha$, $\pr$, and $\rho$ over $V$.
 For $Z$ equal either to $Y$, or $Y_0$, or $V$ we let 
$\hZ=Z\times_{V}U$ (with $U \to V$ being the restriction of $\phi$).
Then  $\alpha_V$, $\pr_V$, and $\rho_V$ induce the morphisms
$\hat \alpha_V : \hY\to \hY_0$, $\hat \pr_V  : \hY_0\to U$, and  $\hat \rho_V : \hY\to U$. 
\enota
\blem\label{main.lb} Let $E$ be an irreducible component of the exceptional divisor of $\alpha_V$,
$C$ be a hypersurface in $E$ such that $\rho_V|_{C}: C \to F^*$ is an unramified covering,
$O$ be an affine neighborhood of $C$ in $Y$, and $f$ be a regular function on $O$
such that $O\cap E=f^*(0)$. Let $\kappa : \hY\to Y$ be the natural projection, $\hE=\kappa^{-1} (E)$,
$\hO=\kappa^{-1} (O)$, and $\hf =\kappa^*f$. Then $\hO\cap \hE=\hf^*(0)$.
Furthermore, for $b \in B$ let  $\hG_b=(b \times \PP^1)\cap U$, $G_b=\phi(\hG_b)$,  $Y^b=\rho_V^{-1} (G_b)$,
and $\hY^b=\hat \rho_V^{-1} (\hG_b)$. Then shrinking $U$ one can suppose 
that $\kappa|_{\hY^b}: \hY^b \to Y^b$ is an isomorphism, $\hY^b$ is smooth, and $Y^b$
meets $E$ transversely in $Y$.
\elem
\bproof
By Lemma \ref{main.l3}
$G_b$ is smooth at $\fo$ and the morphism $\phi_b=\phi|_{\hG_b} : \hG_b \to G_b$ is birational and  \'etale over $\fo$.
Shrinking $U$  we make $G_b$ smooth and $\phi_b$ an isomorphism. 
Let $Y_0^b=\pr_V^{-1} (G_b)$ and $C_0=\alpha_V (C)\subset F_0=\pr_V^{-1} (F^*)$. Since 
$G_b\times \PP^1\simeq Y_0^b \subset Y_0 \simeq V\times \PP^1$ we see
that $Y_0^b$ is smooth and transversal
to $C_0$   by  Lemma \ref{main.l3}.
Shrinking $V$ we can suppose that $V$ does not meet the closure of $T\setminus F$ in $S$.
Thus, $\alpha_V$ is invertible over $Y_0\setminus F_0$.
By Lemma \ref{main.la} 
$Y^b$ is smooth over $\fo$ and it is transversal to $E$. Since  the open
subsets of $Y^b$ and $Y_0^b$ over $V\setminus F^*$ are isomorphic, $Y^b$ is smooth.
Since $\hY^b=\kappa^{-1} (Y^b)$, it is isomorphic to $\hG_b\times_{G_b} Y^b$
and $\kappa|_{\hY^b }: \hY^b \to Y^b$ is an isomorphism.
Therefore, $\hY^b$ is smooth. The transversality of $Y^b$ and $E$ implies that
 $(f|_{Y^b})^*(0)=Y^b\cap O \cap E$. 
 Hence, $(\hf|_{\hY^b})^*(0)=\hY^b\cap \hO \cap \hE$ because $(\kappa|_{\hY^b })^*$ is an isomorphism.
Since this is true for every $b \in B$ we have $\hO\cap \hE=\hf^*(0)$
 which concludes the proof.
\eproof
\blem\label{main.lc}
There is an embedding $\nu : \hY \hookrightarrow  B\times \hY^{b_0}$
 over $B\times \PP^1$
where $ \hY^{b_0}= \hat \rho_V^{-1} (U\cap  (b_0\times \PP^1))$  is a smooth surface.
\elem
\bproof 

Let $\alpha_i: Y_i\to Y_{i-1}$ and  $E_i \simeq C_{i-1}\times \PP^1, \, i=1,2, \ldots $ be as in Lemma \ref{main.la}
for $\alpha=\alpha_V$ and let $E_0=\pr_V^{-1}(F^*)\simeq F^*\times \PP^1$.   
By Lemma \ref{main.l6} the statement is true for $\hY_0$. 
Assume that it is true for $\hY_{m}$ and, thus, we have an embedding 
$\nu_{m} : \hY_{m} \hookrightarrow  B\times \hY_{m}^{b_0}$ over $B\times \PP^1$. 
Let $C_{m}$ be contained in the proper transform of $E_s$ in $Y_{m}$
which by abuse of notation will be also denoted by $E_s$.
Shrinking $F^*$ we can suppose that 
there is an affine neighborhood $O$
of $C_{m}$ in $Y_{m}$ for which
the defining ideal  $I$ of $C_{m}$ is generated by  functions $f$ and $h$ regular on $O$
such that $O\cap E_{s}=f^*(0)$ and the restriction of  $h$ to $E_{s}\simeq C_{s-1} \times \PP^1$ 
is of the form $y^{l}+ g_{l-1}y^{l-1} + \ldots + g_0$ where $y$ is a local coordinate on $\PP^1$ and each
$g_i$ is a regular function on $C_{s-1}$.
Consider the extension of $I$ to the coherent sheaf $\cI$
on $Y_{m}$ whose restriction to $Y_{m} \setminus C_{m}$ is the structure sheaf.
By Lemma \ref{main.la} (iii) $Y_{m+1}$ is the result $\Bl_\cI Y_{m}$ of blowing-up of $Y_{m}$ with respect
to $\cI$. Let $\hat \cI=\phi^{-1}\cI\cdot \cO_{\hY_{m}}$ and, thus, $ \hY_{m+1}\simeq \Bl_{\hat \cI} \hY_{m}$
by \cite[Chapter II, Corollary 7.15]{Har04}.
Consider the natural projection  $\eta : B\times \hY_{m}^{b_0} \to \hY_{m}^{b_0}$ and
the natural embedding 
$\iota: \hY_{m}^{b_0}\simeq b_0\times \hY_{m}^{b_0}\hookrightarrow B\times \hY_{m}^{b_0}$.

{\em Claim}.  There is a coherent sheaf of ideals $\cJ$ on $B\times \hY_{m}^{b_0}$ such that
(a) $\hat \cI=\nu_{m}^{-1} \cJ \cdot  \cO_{\hY_{m}}$ and
(b) $\cJ=\eta^{-1}\cJ_0\cdot \cO_{B\times \hY_{m}^{b_0}}$ where $\cJ_0=\iota^{-1}\cJ \cdot \cO_{ \hY_{m}^{b_0}}$.

Indeed, let $\kappa : \hY_{m} \to Y_{m}$ 
be the natural projection and $L =(\pr_V\circ \alpha_1 \circ \ldots \circ \alpha_m)^{-1} (\fo)\cap E_{s}$.
By Lemma \ref{main.la}(ii)   $L$
is a disjoint union of $n_{s-1}$ lines $L^j\simeq \PP^1$ where $n_{s-1}$ is cardinality of the set
$\{ \fo_{1}, \, \ldots, \fo_{n_{s-1}}\}$  of points   in  $C_{s-1}$ above $\fo$.    
Since $\phi^{-1} (F^*)$ coincides with $\phi^{-1} (\fo)=R\simeq B$ (see Notation \ref{main.n1})
 $\hE_s:=\kappa^{-1}(E_s)\simeq B\times L$ and the restriction of $\hh=\kappa^*h$
to  $B\times L^j$ is the polynomial $y_j^{l}+ d_{j,l-1} y_j^{l-1} + \ldots + d_{j,0}$ 
 with constant coefficients (where $y_j=\kappa_m^*y$ is a local coordinate on $L^j$ and $d_{j,i}=g_i (\fo_j)$).
Consider the ideal $\hI$ 
generated by  $\hh$ and $\hf=\kappa^* f$ on  $\hO =\kappa^{-1}(O)$.
By Lemma \ref{main.lb} $\hO\cap \hE_s=\hf^*(0)$. 
Let us identify $\hY_{m}$ with $\nu_{m}(\hY)$. Then
for a function $\hf^0$ on  $(b_0\times \hY_{m}^{b_0})\cap \hO$ with simple zeros
on $(b_0\times \hY_{m}^{b_0})\cap \hE_s$ the function $\eta^*\hf^0$ has simple zeros on $\hE_s$.
Thus, we can suppose that $\hf=\eta^*\hf^0$. 
Furthermore, let $\hh^{0}=\iota^*\hh$. Then the description of $\hh|_{B\times L^j}$ implies that
 $\hh-\eta^*\hh^0$ is divisible by $\hf$. Hence,
 we can replace $\hh$ by $\eta^*\hh^{0}$ as a generator of $\hI$.
Extend the ideal $\hI$ on $\nu_{m}(O)$  to the coherent  sheaf of ideals $\cJ$  on $B\times \hY_{m}^{b_0}$ 
whose restriction to  $(B\times \hY_{m}^{b_0})\setminus \nu_{m} (O)$ is the structure sheaf.
Let $\cJ_0=\iota^{-1}\cJ \cdot \cO_{ \hY_{m}^{b_0}}$.  Since $\eta^*\hh^0$ and $\eta^*\hf^0$
are generators of
 $\hI$ we see that $\cJ=\eta^{-1}\cJ_0\cdot \cO_{B\times \hY_{m}^{b_0}}$ and
by construction  $\hat \cI=\nu_{m}^{-1} \cJ \cdot  \cO_{\hY_{m}}$ which concludes the proof of the claim.

By  \cite[Chapter II, Corollary 7.15]{Har04}  condition (a) yields the embedding 
$\Bl_{\hat \cI} \hY_{m}\hookrightarrow \Bl_{ \cJ} (B\times \hY_{m}^{b_0})$
and condition (b) yields the isomorphism 
$ \Bl_{ \cJ} (B\times \hY_{m}^{b_0})\simeq B \times \Bl_{\cJ_0}  \hY_{m}^{b_0}$
(with both morphisms over $B\times \PP^1$).
Hence, we have the desired embedding $\nu_{m+1} : \hY_{m+1} \hookrightarrow  B\times \hY_{m+1}^{b_0}$.
By Lemma \ref{main.lb} $ \hY_{m+1}$ is smooth. Now we are done by induction.
\eproof
 \brem\label{main.r1} 
 Let  $\kappa_V: \hY \to Y$ and $\psi_0 : B\times \rho^{-1} (\fo)\to b_0\times \rho^{-1} (\fo)$
be the natural projections.
 Since $\hat \rho_V^{-1} (R)=R \times_{\fo, \phi, \rho} \rho^{-1} (\fo)\simeq  R\times \rho^{-1} (\fo)$
 we have an isomorphism
 $\beta : \hat \rho_V^{-1} (R)\to B\times \rho^{-1} (\fo)$ (because $R=B\times 0$ by Notation \ref{main.n1})
 such that $\beta^{-1} (b_0\times \rho^{-1} (\fo)) =\hY^{b_0} \cap \hat \rho_V^{-1} (R)=:\tE$
 and $\kappa_V|_{\hat \rho_V^{-1} (R)}=\kappa_V\circ \psi$ where 
 $\psi =\beta^{-1} \circ \psi_0 \circ \beta: \hat \rho_V^{-1} (R)\to \tE$.
 The definition of $\eta$ from Lemma \ref{main.lc} implies the commutative diagram
 
 \vspace{0.5cm}
 
  \begin{tikzpicture}\label{it.eq5}   
  \setlength{\unitlength}{1cm}
  
\put(2,0){$ B\times \rho^{-1} (\fo)$}
\put(4.1,0.1){\vector(1,0){1}} 
\put(4.4,0.2){$\simeq$}
\put(5.2,0){$  \hat \rho_V^{-1} (R)$}
\put(6.5,0.1){\vector(1,0){1}} 
\put(6.85,0.2){$\nu$}
\put(7.6,0){$ \nu( \hat \rho_V^{-1} (R))$}
\put(9.5,0.1){\vector(1,0){1}} 
\put(9.8,0.2){$\simeq$}
\put(10.6,0){$ B\times \tE $} 

\put(2.9,-0.2){\vector(0,-1){1}} 
\put(5.8,-0.2){\vector(0,-1){1}} 
\put(8.4,-0.2){\vector(0,-1){1}} 
\put(11.2,-0.2){\vector(0,-1){1}} 
\put(5.9,-0.7){$\psi $} 
\put(8.5,-0.7){$\eta $}

\put(2.4,-1.7){ $\rho^{-1} (\fo)$}
\put(4.1,-1.6){\vector(1,0){1}} 
\put(4.4,-1.5){$\simeq$}
\put(5.6,-1.7){$\tilde E$}
\put(6.5,-1.6){\vector(1,0){1}} 
\put(6.85,-1.5){$\simeq$}
\put(8,-1.7){$ \nu( \tE)$}
\put(9.5,-1.6){\vector(1,0){1}} 
\put(9.8,-1.5){$\simeq$}
\put(11,-1.7){$ \tE $} 
 
     \end{tikzpicture}
  
  \vspace{2cm}
 
 Hence, we have
 \be\label{main.eq1} \kappa_V\circ \nu^{-1}|_{\nu(\hat \rho_V^{-1} (R))}=\kappa_V\circ \nu^{-1} \circ \eta|_{\nu(\hat \rho_V^{-1} (R))}.\ee
\erem 
 The morphism $\alpha : X \to Q$ is invertible over $S\setminus T$ where $T$
is as in Lemma \ref{main.l3}. 
Since $\phi^{-1} (V\setminus T)\supset U\setminus R$ we see that $\hat \alpha_V$
is invertible over $U\setminus R$ which yields an isomorphism 
$\hat \rho_V^{-1} (U\setminus R) \simeq \hat \pr^{-1} (U\setminus R)$ over $U\setminus R$. 

\bprop\label{main.p1}  Let $Z$ be the variety obtained by  the natural gluing of
 $\hY=\hat \rho_V^{-1} (U)$ and $\hQ\setminus \hat \pr^{-1} (R)$ 
 via the isomorphism  $\hat \rho_V^{-1} (U\setminus R) \simeq \hat \pr^{-1} (U\setminus R)$.
Then there is a natural projection $\eta : Z \to B\times \PP^1$ such that  
$Z$ is  isomorphic to $B\times Z_0$ over $Q^{b_0}\times \PP^1$
where  $Z_0=\eta^{-1} (b_0\times \PP^1)$ is a smooth projective surface.
\eprop

\bproof 
Since $\hat \rho_V : \hY \to U$ and $\hat \pr: \hQ\setminus \hat \pr^{-1} (R) \to (B\times \PP^1)\setminus R$
agree over $U\setminus R$
we have $\eta : Z\to B\times \PP^1$.
The  maps from  Lemmas \ref{main.l6} and \ref{main.lc} 
yield the following commutative diagram 
\[\begin{array}{ccl} \hY &  \stackrel{{\nu}}{ \hookrightarrow} & B\times \hY^{b_0}\\
\, \, \, \,\, \downarrow^{\hat \alpha_V}   & 
& \,\,\,\,\,\,  \downarrow^{(\id, \hat \alpha_V^{b_0})}\\
\hY_0 &  \stackrel{{\nu_0}}{ \hookrightarrow} &   B\times \hY^{b_0}_0
\end{array} \]
where $\hat \alpha_V^{b_0} :  \hY^{b_0}\to  \hY^{b_0}_0$ is the restriction of $\hat \alpha_V$,
and the vertical arrows are biregular over $U\setminus R$.
By Lemma \ref{main.l5} there is also an isomorphism of $\hQ\setminus \hat \pr^{-1} (R)$ and $B\times (\hQ^{b_0} \setminus \hat \pr^{-1} (R))$
which agrees with $\nu_0$ over $U\setminus R$.
Since $Z_0$ is the result of gluing of $\hQ^{b_0} \setminus \hat \pr^{-1} (R)$ and $\hY^{b_0}$ over $U\setminus R$
we have the desired conclusion.
\eproof
\blem\label{main.l8} Let the notations of Proposition \ref{main.p1} hold,
  $\hX=(B\times \PP^1)\times_{S,\phi,\rho} X$, 
$\hD=(B\times \PP^1)\times_{S,\phi,\rho}D$ (where $D$ is as in Lemma \ref{main.l2}), 
$\hat \rho : \hX \to B\times \PP^1$ and  $\kappa : \hX \to X$ be the natural projections,
and $\hv=\kappa^{-1} (v)$
(since $\phi$ is invertible over
a neighborhood of $s_0$ we see that $\kappa$ is invertible
over $v$ and, thus, $\hv$ is a singleton).
Then there is a morphism $\gamma : \hX \to Z$ for which

{\rm (i)} $M=\gamma (\hD \setminus \hat \rho^{-1} (R))$ is   a codimension at  least 2 subvariety in $Z$
such that for every $b \in B$ the set $M_b=Z^b \cap M$ is finite where $Z^b =\eta^{-1} (b\times \PP^1)$;

{\rm (ii)} $\gamma$ is biregular over $Z\setminus M$ (in particular,
$\gamma$ is biregular over a neighborhood of $\gamma (\hv)$).

\elem
\bproof  Let $\hat \alpha : \hX \to \hQ$ be the morphism  induced by $\alpha$.
Consider the restriction $\gamma_1$ of $\hat \alpha$ over $(B\times \PP^1) \setminus R$
and the morphism $\gamma_2: \hat \rho^{-1}(U) \to \hat \rho_V^{-1} (U)$ induced by 
the isomorphism $ \rho^{-1}(U) \to \rho_V^{-1} (U)$ (that is, $\gamma_2$ is the identity isomorphism).  
Since $\alpha^{-1} (S\setminus T)$ is naturally isomorphic to $\pr^{-1} (S\setminus T)$ 
 we see that $\gamma_1$ and $\gamma_2$ agree over $U\setminus R$
which yields a morphism $\gamma : \hX  \to Z$ which is biregular over $U$.
By Lemma \ref{main.l3} $\phi$ is invertible over a neighborhood $W$ of $s_0$ in $S$
and by  Lemma \ref{main.l2} (ii) $\rho^{-1} (W)$ does not meet $D$.
Hence, $\hD$ does not meet $\hat \rho^{-1} (\phi^{-1} (W))$.
Shrinking $B$ we can suppose that  $\phi (b\times \PP^1)$ meets $W$ 
for every $b \in B$. Then the
set $M_b$ is contained in the union of a finite number of
fibers of  $\eta|_{Z^b} : Z^b\to \PP^1$.
Furthermore, by Lemma \ref{main.l2} (i) $M_b$ cannot contain a fiber of $\eta_b$.
Hence, we have statement (i). 
Since  $\hat \alpha$ and $\gamma$ coincide over $(B\times \PP^1) \setminus R$
statement (i) implies that  $\gamma$ is biregular on the complement of
$\hD \setminus \hat \rho^{-1} (R)$ in $\hX$.
Since $v \notin D$ we have $\hv \notin \hD$ which concludes (ii) and the proof.
\eproof
\bprop\label{main.p2}  Let $\cD$ be a nontrivial line bundle on a smooth complete rational
variety $X$ and $v\in X$ be a point that possesses  a neighborhood isomorphic to
an open subset of $\A^n$. Let $\pi : \tX \to X$ be a morphism from
a smooth rational variety $\tX$ locally invertible over a neighborhood of $v$, 
$\tilde \cD=\pi^*\cD$, and $\tv=\pi^{-1} (v)$.
Suppose that  for every triple $(X, \cD, v)$ as above the triple $(\tX, \tilde \cD, \tv)$ can be chosen so that $\cP$-property
is valid for it whenever $\dim X =n-1$.
Then the same remains true if $\dim X=n$.
\eprop 

\bproof  Since $v$ has a neighborhood isomorphic to $\A^n$ we can suppose that
the conclusions of Lemma \ref{main.l2} hold. 
Let $F$ be an irreducible component of $\rho (\supp D')$.
Consider the triple $(S,\cF, s_0)$ where $\cF=\cO_S(F)$.
By the assumption there exists a morphism $\pi : \tS \to S$ invertible over a neighborhood of $s_0$
such that $\cP$-property holds for the triple $(\tS, \tilde \cF, \ts_0)$ where $\tilde \cF =\pi^* \cF$ and
$\ts_0=\pi^{-1} (s_0)$. Let $\tX=X\times_S \tS$, $\tilde \cD$ be the lift of $\cD$ to $\tX$,
and $\tv$ be the preimage of $v$ in $\tX$. It suffices to prove
$\cP$-property for $(\tX, \tilde \cD, \tv)$. Changing notations we can suppose
that  $(X,  \cD, v)=(\tX, \tilde \cD, \tv)$ and $(S,\cF, s_0)=(\tS, \tilde \cF, \ts_0)$. 
Thus, we have $\phi : B\times \PP^1 \to S$ and $\fo\in F$
as in Lemma \ref{main.l3}. For $b_0\in B$
such that $s_0\in \phi (b_0\times \PP^1)$ one has $v\in \rho^{-1} ( \phi (b_0\times \PP^1))$. 
Let the notations of Proposition \ref{main.p1} and Lemma  \ref{main.l8} hold.

{\em Claim 1.} There is a morphism $\hat \psi : W\times \PP^1 \to \hX$
invertible over $\hv$ and an irreducible component $\hH_0$ 
 of $\hD'=\kappa^{-1} (\supp D')$
such that every curve $\hat \psi (w\times \PP^1),\, w\in W$  intersects $\hD'$  at one point $w''\in \hH_0$ only and
this intersection is transversal.

Let  $\hH=\hat \rho_V^{-1} (R) \cap \hD'$ and $L=\gamma (\hH)\cap Z_0$.
Note that the collection $\{ Z_0\to Q^{b_0}, \hat \alpha_V(\hv), L\}$ satisfies the assumptions for
the collection $\{ \beta: Y \to P, u_0 , T\}$ in Lemma \ref{it.l5}.
Hence, by Proposition \ref{it.p3} there is an irreducible component $L_0$ of $L$
for which one can find a birational morphism $\psi_0 : B_1\times \PP^1 \to Z_0$ invertible over
a neighborhood of $\gamma (\hv)$ in $Z_0$ such that $\psi_0 (B_1 \times \PP^1)\cap M_{b_0}=\emptyset$,
the set  $\psi_0 (B_1 \times \PP^1) \cap L$  is a general point $\fo_1$ of $L_0$,
and  every curve $\psi_0 (b_1 \times \PP^1), \, b_1\in B_1$ intersects $L_0$  transversely at $\fo_1$.     
 By Proposition \ref{main.p1} we can identify $Z$ with  $B\times Z_0$. 
  Then  the morphism 
  $ \psi: B\times B_1\times \PP^1\to Z, \, (b,b_1,t) \mapsto (b, \psi_0(b_1,t))$ 
 is invertible over a neighborhood of  $\gamma (\hv)$ in $Z$.
 Let $b_1^0\in B_1$ be such that $\gamma (\hv) \in \psi (b_0\times b_1^0 \times \PP^1)$.
 Replacing $B\times B_1$ with a  neighborhood $W$ of $(b_0,b_1^0)$
we can suppose that $\psi (W\times \PP^1)$ does not meet $M$ and
every curve $ \psi (w \times \PP^1), \, w=(b, b_1)\in W $ 
 intersects $\gamma (\hH)$  at the point 
 $w'=(b, \fo_1)\in B\times L_0$ only
 with the intersection being  transversal.
 By Lemma \ref{main.l8} $\gamma (\hD) \subset \gamma (\hH) \cup M$.  Hence,
by Lemma \ref{main.l8} (iii) one has a lift $\hat \psi : W\times \PP^1 \to \hX$ of $\psi$
invertible over $\hv$
such that every curve $\hat \psi (w\times \PP^1)$  intersects $\hD'$  at $\hw''=\gamma^{-1}(w')$ only
 with the intersection being  transversal.  Furthermore,  $\hw''$  is  contained in the irreducible
 component $\hH_0=\gamma^{-1} (B\times L_0)$ of $\hH$ which concludes the proof of Claim 1. 
 
{\em Claim 2.}
 There is a morphism $\chi : W\times \PP^1 \to X$
invertible over $v$, an irreducible component $H$
 of $D'$, and a general point $\breve \fo \in H$
such that every curve $\chi (w\times \PP^1),\, w\in W$  intersects $\supp D'$  at $\breve \fo$  only and
this intersection is transversal.

 Note that $\hH=\kappa^{-1}(E)$ where $E= \supp D'  \cap \rho^{-1} (\fo )$. 
 Hence, $\hH \simeq B\times E$ by Remark \ref{main.r1},  and, consequently,
  $\gamma (\hH) \simeq B\times E$.
 By Formula \eqref{main.eq1} $\kappa\circ \gamma^{-1}$ sends
 every point $(b,e)\in \gamma (\hH)$ to $e\in E$.
In particular,  $\kappa \circ \gamma^{-1} (w')=\kappa \circ \gamma^{-1} (b,\fo_1)=\breve \fo$ 
 is a general point of $E_0=\kappa (\hH_0)$ independent of $b$.
 Let  $H$ be an irreducible component of $D'$
for which $E_0=\kappa (\hH_0)$ is an irreducible component of $H\cap \rho^{-1} (\fo)$. 
 Since $\rho (E_0) =\fo$ is a general point of $F$ we see that $\breve \fo$ is a general point of $H$. 
For  $ \chi=\kappa \circ \hat \psi: W\times \PP^1\to X$ we have  $\chi (W \times \PP^1)\cap \supp D' =\breve \fo$.
Since $\kappa$ is invertible over $v$ (see Lemma \ref{main.l8}) and by Claim 1
 $ \hat \psi$ is invertible over $\hv$ we see that $\chi$ is invertible over $v$.
 It remains to show that every curve $\chi (w \times \PP^1)$  intersects $\supp D'$  at  $\breve \fo$ transversely.
 For $w=(b,b_1) \in W\subset B\times B_1$ the curve $\hat \psi (w \times \PP^1)$ is contained in 
 the surface $\hX^b=\hat \pr^{-1} (b \times \PP^1)$
 and $\chi (w \times \PP^1)$  is contained in the surface $X^b= \pr^{-1} (G_b)$ where $G_b=\phi (b\times \PP^1)$.
 By Lemma \ref{main.lb} $X^b$ meets $H_0$ transversely and $\kappa|_{\hX^b}: \hX^b\to X^b$ is birational morphism
 \'etale over $X^b \cap H_0$. Since $\hat \psi|_{w\times \PP^1} : w\times \PP^1 \to \hX^b$ is the lift of
 $\psi|_{w\times \PP^1} : w\times \PP^1 \to Z^b$ we see
 that $\hat \psi (w\times \PP^1)$ meets $\hX^b \cap \hH_0$ transversely in $\hX^b$. 
Hence, $\chi(w\times \PP^1)$ meets $X^b \cap H_0$ transversely  in $X^b$. Therefore,  $\chi(w\times \PP^1)$ 
 is transversal to $H_0$  in $X$
 which yields Claim 2 and, therefore, concludes the proof.
 \eproof
\proof[Proof of Theorem \ref{main.t1}]\label{main.proof} 
By Example \ref{pre.exa1} we can suppose that $\dim X\geq 2$.
By Proposition \ref{it.p3} and Lemma 
\ref{it.l4} for a smooth projective surface $X$ and a triple $(X,\cD, v)$ one can find a
triple $(\tX, \tilde \cD, \tv)$ as in Proposition \ref{main.p2} which satisfies the $\cP$-property.
Proposition \ref{main.p2} and induction imply that the same is true when $X$ is a complete uniformly
rational variety of any dimension. Thus, we can suppose that  there is
 a birational morphism $\phi : B\times \PP^1 \to X$ invertible over
a neighborhood of $v$ such that every curve $ \phi (b \times \PP^1), \, b \in B$ meets only
one irreducible component of  $D'$.
Hence, the restriction of $\cD$ to any of these curves is a nontrivial line bundle.
Every complete uniformly rational variety is elliptic by Theorem \ref{int.t1}.
Now Corollary \ref{pre.c1}  implies the desired conclusion. \hfill $\square$\\

While discussing  principal $\G_m$-bundles  we suppose a priori that they are locally  trivial in the  \'etale topology
 but first let us state an unpublished observation of L\'arusson.
 
 \bprop\label{main.p3} {\rm (L\'arusson)} Let $q: Y' \to Y$ be a morphism of smooth varieties
 and $\iota : Y \hookrightarrow Y'$ be an embedding such that $q\circ \iota =\iota$.
 If $Y'$ is elliptic, then so is $Y$. 
 \eprop
 
 \bproof
Let $(E',p',s')$ be a dominant spray on $Y'$. Identify $Y$ and $\iota(Y)$ and consider
the restriction   $p: E \to Y$  of $p'$. Let $s=q\circ s'|_E: E\to Y$. Then $(E,p, s)$ is a dominant spray on $Y$.
\eproof

\bcor\label{main.c1} Let $X$ be a complete stably uniformly rational variety and $\beta : Y\to X$
be a nontrivial principal $\G_m$-bundle. Then $Y$ is elliptic.
\ecor

\bproof  Suppose first that $X$ is rational.
The morphism $\beta$ extends to a proper morphism $\bY \to X$ which is a $\PP^1$-bundle
such that the divisor $G=\bY \setminus Y$ is an unramified two-sheeted covering of $X$.
Note that $G$ cannot be irreducible since  $X$ is simply connected by \cite[Lemma 2]{Ser59}.
Hence, $G$ consists of two sections. 
Thus, $Y$ can be viewed as the complement to a zero
section of a nontrivial line bundle over $X$. By Theorem \ref{main.t1} this implies ellipticity of $Y$. 
In the general case consider the natural embedding of an elliptic
$X'=X\times \A^m$ into $X''=X\times \PP^m$.
The Gromov Localization Lemma (\cite[82, 3.5.B, 3.5.C]{Gro89} and \cite[Theorem 6.2]{For23},
see also \cite[Corollary 2.2. and Remark 2.3]{KZ24} for a correction of Gromov's proof)
implies that $X''$ is also elliptic. Note that $\beta$ induces
a $\G_m$-bunlde  $\beta'' : Y''\to X''$.  As before we see that $Y''$ can be treated
as the complement to the zero section of a nontrivial line bundle on $X''$. 
Thus,  $Y''$ is elliptic. 
 Since $Y$ is a retract of $Y''$, by Proposition \ref{main.p3}  $Y$ is elliptic 
and we are done.
\eproof
M. Zaidenberg suggested the following strengthened version of Question 1 in Introducation
to which the author does not know the answer.\\

{\bf Question 2} (Zaidenberg). Let $X$ be a smooth complete elliptic variety that
verifies the strengthened curve-orbit property and $\cD \to X$ be a nontrivial
line bundle with the zero section $Z_\cD$. Is $Y=\cD \setminus Z_\cD$ elliptic? \\

It is worth mentioning that there are non-rational unirational varieties that
satisfy the strengthened curve-orbit property (see \cite{KZ24}).\\

{\em Acknowledgements.} The author is grateful to M. Zaidenberg for useful consultations
and the referee for very valuable suggestions.


%

\end{document}